\documentclass[fleqn,10pt]{llncs}

\usepackage{algorithm}
\usepackage{algpseudocode}
\usepackage{amsmath}
\usepackage{amssymb}
\usepackage{color}
\usepackage{graphicx}
\usepackage[utf8]{inputenc}
\usepackage{subfig} 
\newcommand{\nats}{\mathbb{N}}

\newcommand{\real}{\mathbb{R}}
\newcommand{\A}{B}

\newcommand{\din}{k}
\newcommand{\dout}{\hat{k}}
\newcommand{\msp}{~\vert ~}
\newcommand{\pa}{\pi^{pa}_t}
\newcommand{\ra}{\pi^{ra}_t}
\newcommand{\argmax}{\mathop{\mathrm{argmax}}}

\providecommand{\floor}[1]{\left \lfloor #1 \right \rfloor }

\newcommand{\Sref}[1]{Section~\ref{#1}} 
\newcommand{\Fref}[1]{Figure~\ref{#1}} 
\newcommand{\Eref}[1]{Equation~\ref{#1}}

\title{Estimating Formation Mechanisms and Degree Distributions in
  Mixed Attachment Networks}
\author{Jan A. Medina, Jorge
  Finke, Camilo Rocha}
\institute{Department of Electrical
    Engineering and Computer Science\\Pontificia Universidad
    Javeriana\\Cali, Colombia.}

\begin{document}

\maketitle

\begin{abstract}
Our work introduces an approach for estimating the contribution of
attachment mechanisms to the formation of growing networks. We present
a generic model in which growth is driven by the continuous attachment
of new nodes according to random and preferential linkage with a fixed
probability. Past approaches apply likelihood analysis to estimate the
probability of occurrence of each mechanism at a particular network
instance, exploiting the concavity of the likelihood function at each
point in time. However, the probability of connecting to existing
nodes, and consequently the likelihood function itself, varies as
networks grow. We establish conditions under which applying likelihood
analysis guarantees the existence of a local maximum of the
time-varying likelihood function and prove that an expectation
maximization algorithm provides a convergent estimate. Furthermore,
the in-degree distributions of the nodes in the growing networks are
analytically characterized. Simulations show that, under the proposed
conditions, expectation maximization and maximum-likelihood accurately
estimate the actual contribution of each mechanism, and in-degree
distributions converge to stationary distributions.
\end{abstract}

\noindent{\it Keywords\/}: Complex networks, Network model, Statistical inference.

\section{Introduction}

The aim of a wide range of network models is to provide a framework to
understand how linkage mechanisms for establishing links give rise to
particular topological properties, including degree
distributions~\cite{Barabasi1999,Clauset2009},
clustering~\cite{Szabo2004}, average path lengths~\cite{Albert2002},
and community partitions~\cite{Newman2006}. There has been a
continuous and significant effort directed at formalizing such
mechanisms and their role in the formation of networks like the
Internet, the world wide web, and co-authorship
associations~\cite{Sheridan2008,Bianconi2000,Shao2006,Callaway2001,Amaral2000,Krapivsky2000,Dorogovtsev2000,Liu2002,Dorogovtsev2000structure,Ke2004,Jackson2007,Barabasi1999,Tadic2001}.

The work in~\cite{Barabasi1999} explains the emergence of power law
degree distributions as an outcome of two mechanisms, the addition of
new nodes to the network (a growth mechanism) and the preference of
new nodes to connect to nodes with a high degree (an attachment
mechanism). The two mechanisms combined yield a degree distribution
that asymptotically follow a power law. If the growth mechanism is
combined with a uniform random attachment mechanism, with no
preference for nodes with a high degree, then the resulting network
follows an exponential function~\cite{Shao2006}. However, a recent
study shows that from nearly $1000$ empirical networks less than 5\%
exhibit pure power law or exponential relationships~\cite{Broido2018}.

Empirical networks have in general a richer diversity in structure,
which suggests that, as networks grow, multiple mechanisms of
attachment contribute to the resulting degree
distributions~\cite{Wang2012,Zhang2015}. Citation networks, for
example, have degree distributions which obey neither exponential nor
power law, but mixed distributions ~\cite{Redner1998}. Evidence for
such degree distributions are also found in opinion
networks~\cite{Said2017}, and protein-protein interaction
networks~\cite{Giot2003,Li2004}.

It has been observed that networks in diverse domains (e.g.,
biological processes and social interaction) can have similar
underlying generative mechanisms that make them hard to distinguish
based purely on network structure~\cite{Kansuke2017}. Therefore, there
has been especial interest in finding connections between the
structure of networks and their underlying formation
mechanisms~\cite{Strogatz2001,Newman2003}. A number of models have
been introduced to explain how attachment mechanisms, in particular,
give rise to distinct topological properties.  An important step in
developing these models entails the problem of assessing the
plausibility of each mechanism. Previous research tries to estimate
the contribution of a mechanism using maximum-likelihood methods. Such
methods determine the optimal estimate that best describes the
contribution of each mechanism based on the number of new edges
established at some point in time~\cite{Wang2012,Zhang2015}.  The
work in~\cite{Zhang2015} assigns an adjustable weight to each
mechanism. The optimal estimate represents the set of weights that
maximize the likelihood of all new edges. However, as the network
grows, applying standard maximum-likelihood estimation does not
produce a consistent estimate over time.

Our work focuses on understanding how to obtain a convergent estimate
of the contribution of multiple mechanisms that influence the
evolution of growing networks.  In particular, it aims to provide an
analytical framework to evaluate the contribution of two attachment
mechanisms, namely, random and preferential attachment. We extend the
method in~\cite{Wang2012} and~\cite{Zhang2015} by presenting
conditions to guarantee the existence of a realistic
maximum-likelihood estimate. Moreover, an expectation maximization
algorithm is applied to evaluate the contribution of these two
attachment mechanisms in one simulated network and one empirical
citation network.

The contributions of this paper are the following. First, we
characterize the roots of the likelihood function regarding the
network parameters (Lemmata~\ref{lem1} and \ref{lem2}). Second, we
present conditions under which the likelihood function has a maximum
(Theorem~\ref{thm1}) and an algorithm to estimate the contribution of
the two attachment mechanisms (preferential and random
attachment). Third, we use a discrete-time approach to characterize
the in-degree distribution as a function of the parameters and
contribution of each mechanisms (Theorem~\ref{thm2}).  Fourth, we show
that the dynamics of the in-degree distribution converges to a
stationary distribution
(Corollary~\ref{cor.indegree.behavior}). Finally, we verify that the
estimate of the contribution of random and preferential attachment
yields a theoretical in-degree distribution that resembles empirical
distributions of citation networks.

The remainder of the paper is organized as follows. \Sref{section2}
introduces the mixed attachment model. \Sref{sec.maxlike} presents the
proposed estimation approach. \Sref{sec.em} overviews the estimation
based on expectation-maximization algorithm. \Sref{section4}
characterizes the in-degree distribution of the model. \Sref{section5}
presents simulation results. \Sref{section6} draws the concluding
remarks and some future work.

\section{The Network Model}
\label{section2}

The network model used in this paper is an
extension of the network model in~\cite{Shao2006}, which supports directed 
networks and includes a response mechanism. It consists of three main
mechanisms, namely, growth, attachment, and response. By
\textit{growth} we mean that the number of nodes in the network
increases by one at each time step. \textit{Attachment} refers to the
fact that new nodes tend to connect to existing nodes, while
\textit{response} refers to the fact that existing nodes tend to
connect to new nodes. In mathematical terms, the network model is
parametric in a probability $\alpha$ and natural numbers $m$ and $\hat{m}$ 
governing the attachment and response mechanisms. Internally, the attachment mechanism
creates $m>0$ outgoing edges from the new node and is characterized as
a Bernoulli trial with parameter $\alpha$, where $\alpha$ represents
the probability of establishing a new edge by preferential attachment
and $1 - \alpha$ by random attachment. The response mechanism creates
$\hat{m}\geq 0$ incoming edges from the existing nodes to the new node
by random attachment.

A network is represented as a directed graph $G_t=(V_t,E_t)$ with
nodes $V_t$ and edges $E_t \subseteq V_t \times V_t$. A pair $(u,v)\in
E_t$ represents a directed edge from a \textit{source} node $u$ to a
\textit{target} node $v$. The expressions $\din_t(u)$ and $\dout_t(u)$
denote, respectively, the \textit{in-} and \textit{out-degree} of node
$u\in V_t$. Moreover, $n_t$ and $e_t$ denote the \textit{number of
  nodes} and the \textit{number of edges} in the network at time $t$,
respectively (i.e., $n_t=|V_t|$ and $e_t=|E_t|$).

\begin{definition}\label{def.netmodel.netmodel}
  The \textit{algorithm used in the network model} goes as follows:

  \begin{enumerate}
  \item Growth: starting from a seed network $G_0$, at each time step
    $t>0$, a new node is added with $m$ outgoing edges that link the
    new node to $m$ different nodes already present in the network and
    $\hat{m}$ incoming edges that link $\hat{m}$ different nodes
    already present in the network to the new node.
  \item Attachment: when choosing the $m$ nodes to which the new node
    connects, we assume that the probability $\pi_t(v\mid \alpha)$
    that the new node will be connected to node~$v$ is given by
\begin{equation}
  \pi_t(v\mid \alpha)= \alpha\pa (v) + (1-\alpha)\ra (v) \label{eq1}
\end{equation}
where the probability of establishing an outgoing edge from the new
node to the existing node $v$ due to preferential attachment is given
by
  \begin{equation}
    \pa(v)=\frac{\din_{t-1}(v)}{e_{t-1}} \label{eq2}
  \end{equation}
and due to random attachment by
  \begin{equation}
    \ra(v)=\frac{1}{n_{t-1}}. \label{eq3}
  \end{equation}
\item Response: when choosing the $\hat{m}$ nodes that connect to the
  new node, we assume that the probability $\eta_t(v)$ that an
  existing node $v$ connects to the new node is given by
\begin{equation}
  \eta_t(v)= \frac{1}{n_{t-1}}. \label{eq3.1}
\end{equation}
  \end{enumerate}
\end{definition}

For the growth and response processes to be well-defined, the
algorithm in Definition~\ref{def.netmodel.netmodel} assumes that the
seed network has at least $\max\{m,\hat{m}\}$ nodes. Although it is
not required by the algorithm, it will be further assumed that the
seed network does not have self-loops and each node has at least one
incoming edge. In this way, the networks generated by the algorithm do
not contain self-loops and each node has non-zero in-degree.  As
explained before, the attachment mechanism is characterized as a
Bernoulli trial with parameter $\alpha$, representing the probability
of establishing a new edge via preferential attachment. For
preferential attachment, the probability of establishing a new edge
depends on the in-degree of the target node, meaning that the
probability for an existing node of becoming a target node is directly
proportional to its in-degree. For random attachment, the probability
of establishing an edge to a target node follows a discrete uniform
distribution. The response mechanism selects from the set of existing 
nodes uniformly at random which nodes connect to the new node.

\begin{figure}[htbp]
  \centering
  \includegraphics[width=200pt]{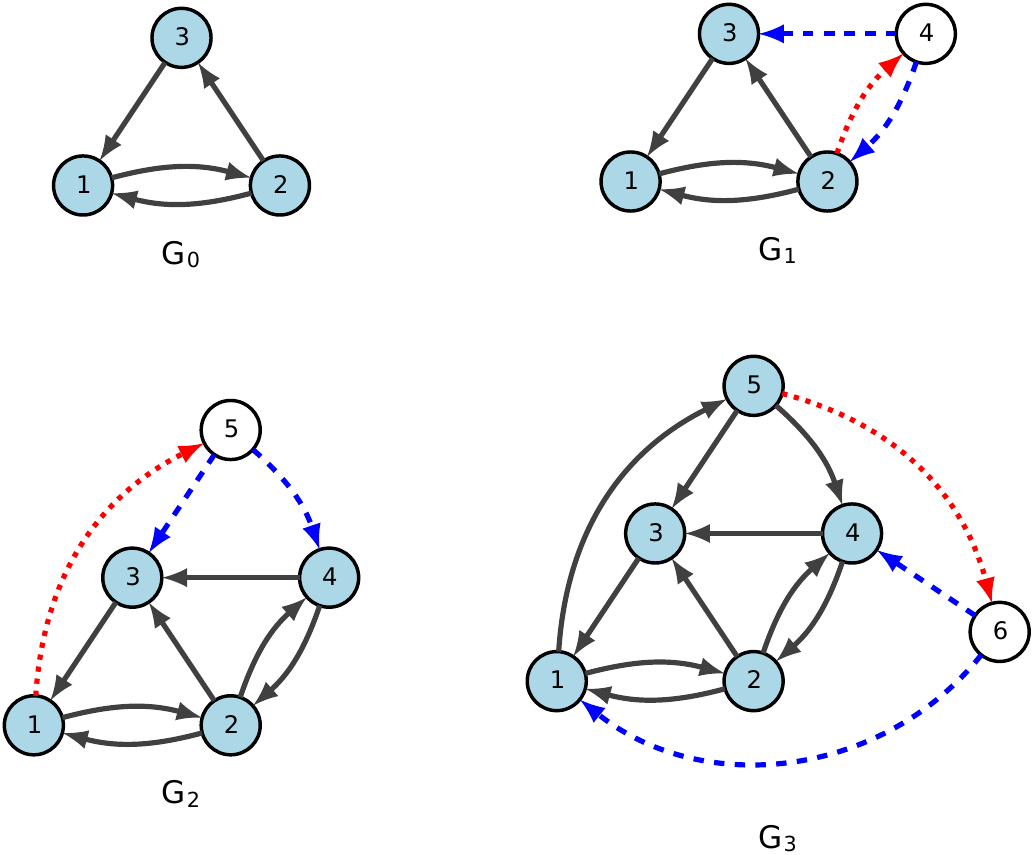}
  \caption{Network evolution from a given seed network $G_0$ in the
    network model. At each time step $t \geq 1$, a new node is added
    with $m=2$ outgoing edges (dashed line depicted in blue) and $\hat{m}=1$
    incoming edges (dotted line depicted in red).}
\label{fig:network}
\end{figure}

\Fref{fig:network} illustrates the evolution of a seed network $G_0 =
(V_0, E_0)$ at time steps $t = 1, 2, 3$ in the network model with
$m=2$ and $\hat{m}=1$ (for the sake of simplicity in the illustration,
the probability $\alpha$ is omitted). The seed network $G_0$ has nodes
$V_0 = \{1, 2, 3\}$ and edges $E_0 = \{(1,2), (2,1), (2,3), (3,1)\}$.
At time $t=1$, the set of nodes grows by adding the new node 4, and by
creating the new edges $(4,2), (4,3), (2, 4)$ ($m=2$ outgoing edges by
the attachment mechanism and $\hat{m}=1$ incoming edges by the
response mechanism). At time $t=2$, the new node $5$ establishes $m=2$
outgoing edges and $\hat{m}=1$ incoming edges, as it is the case for
the new node $6$ at time $t=3$.

\section{Maximum Likelihood Analysis}
\label{sec.maxlike}

Consider the problem of determining the values of the $\alpha$, $m$,
and $\hat{m}$ parameters used in the process of generating a sequence
of netwroks $\left(G_t\right)_{t\leq T}$ by the algorithm in
Definition~\ref{def.netmodel.netmodel} (see \Sref{section2}), for some
seed network $G_0$ and $T\in \nats$.  Inspection on such a sequence
could be enough for establishing the values of $m$ and
$\hat{m}$. However, determining the value of $\alpha$ is, in general,
a tall order. The purpose of this section is two-fold. On the one
hand, it argues that recently published maximum likelihood analysis
approaches (e.g., in~\cite{Zhang2015}) do not produce reasonable
estimates in some cases. On the other hand, it presents a new approach
to provide better estimates and sufficient conditions under which such
an analysis is feasible.

The likelihood of creating an edge from the new node $u$ at time $t$
to a node $v\in V_{t-1}$ is given by~\cite{Zhang2015}:
\begin{eqnarray}
\pi_t(v\mid \alpha) &=& \alpha\pa(v) + (1-\alpha)\ra(v) \nonumber\\
&=& \alpha\frac{\din_{t-1}(v)}{e_{t-1}}+(1-\alpha)\frac{1}{n_{t-1}}\nonumber\\
&=& \alpha\left(\frac{\din_{t-1}(v)}{e_{t-1}}-\frac{1}{n_{t-1}}\right)+\frac{1}{n_{t-1}}.\label{eq4}
\end{eqnarray}
Since preferential attachment depends on the in-degree of the target
node, \Eref{eq4} can be written as a function of the in-degree of node
$v$ with $k = \din_{t-1}(v)$:
\begin{equation}
\pi_t(k\mid \alpha) = \alpha\left(\frac{k}{e_{t-1}}-\frac{1}{n_{t-1}}\right)+\frac{1}{n_{t-1}}.\label{eq5}
\end{equation}

\begin{definition}\label{def.maxlike.indegs}
  For $t \in \nats$, let:
  \begin{itemize}
    \item \emph{$A_t$ be the multiset of in-degrees} of nodes selected
      at time $t$ by the attachment mechanism in
      Definition~\ref{def.netmodel.netmodel} and
    \item $\A_t = \{A_1,A_2,\ldots,A_t\}$ be the collection of all
      multisets of in-degrees of nodes selected up to time $t$.
  \end{itemize}
\end{definition}

At each time $t$, the family of sets $\A_t$ is a random sample of size
$t$. Since the elements of $\A_t$ are independent and identically
distributed (i.i.d.), it follows from~\Eref{eq5} that the likelihood
function can be given by~\Eref{eq6}:
\begin{eqnarray}
f_t(\alpha)&=&\prod\limits_{i=1}^t\prod\limits_{j=1}^{|A_i|}\pi_i(k_{ij}\mid \alpha)\nonumber\\
&=& \prod\limits_{i=1}^t\prod\limits_{j=1}^{|A_i|}\left[\alpha\left(\frac{k_{ij}}{e_{i-1}}-\frac{1}{n_{i-1}}\right)+\frac{1}{n_{i-1}}\right],\label{eq6}
\end{eqnarray}
where $k_{ij}$ denotes the $j$-th element in $A_i$ (without loss of
generality, the multiset $A_i$ is assumed sorted in ascending
order). Note that $f_t(\alpha)$ is a polynomial in the indeterminate
$\alpha$ and has order at most $mt$. If all $k\in A_i$ satisfy $k
n_{i-1} - e_{i-1} \neq 0$ for $1\leq i\leq t$, then the order of
$f_t(\alpha)$, denoted $\deg{f_t}(\alpha)$, is $mt$.

\begin{definition}\label{def.maxlike.likelihood}
For a network $G_t$, the function $f_t$ is called the
\emph{$G_t$-likelihood function}. The \textit{maximum likelihood
  estimator} is defined as
\[\hat{\alpha}=\argmax_{\alpha\in(0,1)}f_t(\alpha).\]
\end{definition}

Consider the complete directed graph $G_0$ with 3 vertices, and a
sequence $(G_t)_{t\leq 2000}$ generated by the algorithm in
Definition~\ref{def.netmodel.netmodel} with parameters $m = 5$,
$\hat{m} = 3$, and $\alpha=0.6$.

\begin{figure}[tbhp]
	\centering
    \includegraphics[width=170pt]{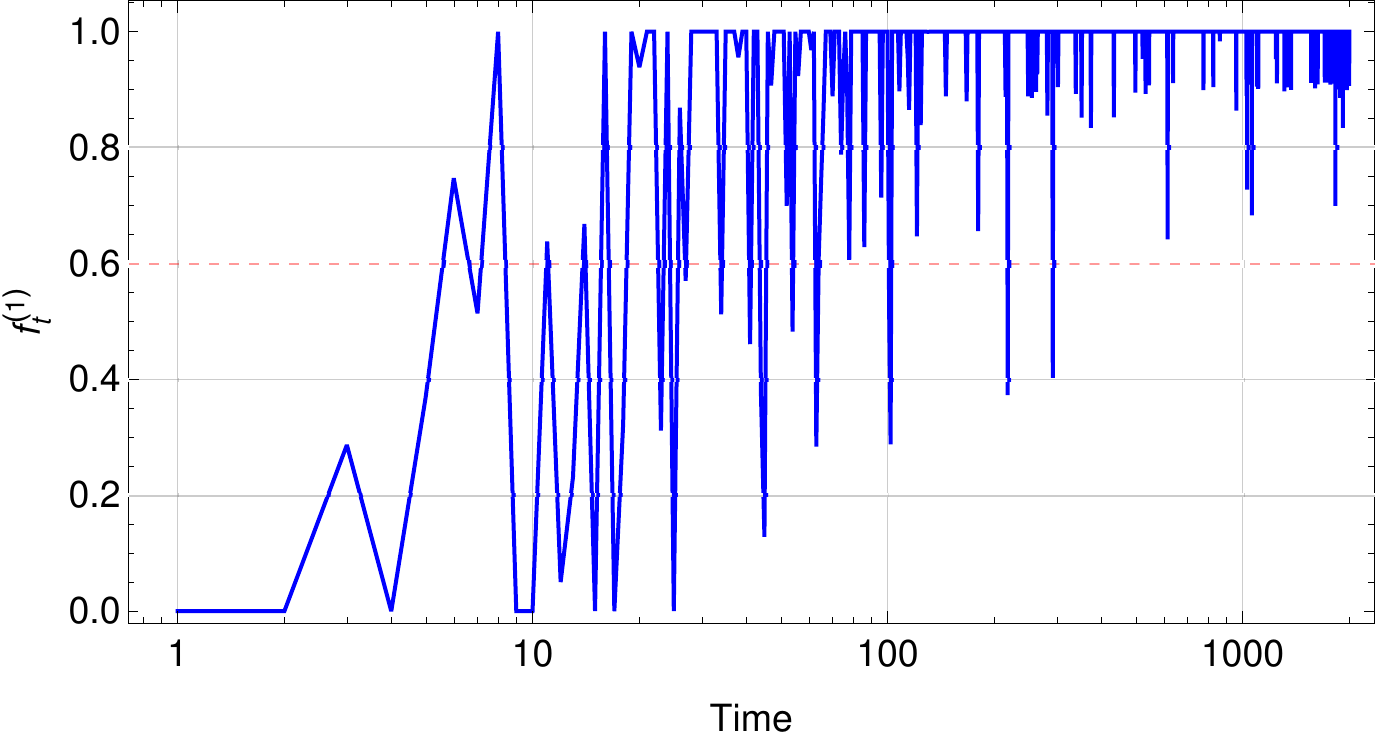}~
    \includegraphics[width=170pt]{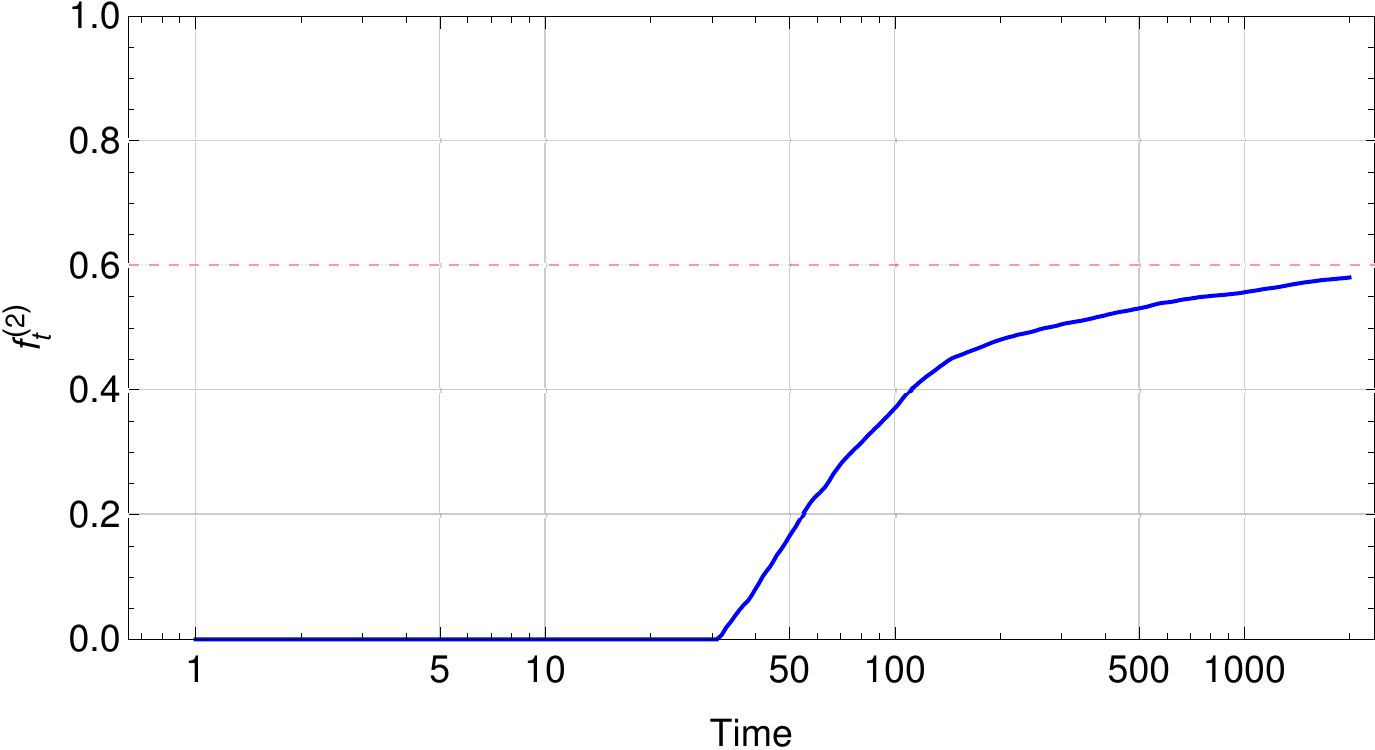}
\caption{Behavior of the parameter estimated by the likelihood functions $f^{(1)}_t$ and $f^{(2)}_t$, 
with $m=5$, $\hat{m} = 3$, and $\alpha = 0.6$.}\label{fig1.maxlike}
\end{figure}

At each time $t$, likelihood analysis is applied to the multiset $A_t$
and to the collection $\A_t$. Let the likelihood functions
$f^{(1)}_t(\alpha)$ and $f^{(2)}_t(\alpha)$ be defined as:
\[f^{(1)}_t(\alpha) = \prod\limits_{k\in A_t}\pi_t(k\mid\alpha)\qquad\text{and}\qquad f^{(2)}_t(\alpha) = \prod\limits_{i=1}^t\prod\limits_{j=1}^{|A_i|}\pi_i(k_{ij}\mid \alpha).\]
At each time, an estimation of $\alpha$ is provided by each likelihood
function. \Fref{fig1.maxlike} shows that the standard
maximum-likelihood estimation with $f^{(1)}_t$ is not capable of
producing a consistent estimate over time (this witnesses a
counter-example to the approach in~\cite{Zhang2015} for likelihood
analysis). However, the proposed approach developed in the rest of
this section, which supports the maximum-likelihood estimation with
$f^{(2)}_t$, yields a better estimate. The remaining of this section
establishes conditions under which the estimator based on the
$G_t$-likelihood function produces a consistent estimate over time for
the network model in Section~\ref{section2}.

The next goal is to identify some key relationships between the number
of nodes, the number of edges, and the number of roots of the
$G_t$-likelihood function $f_t$, and then provide conditions under
which it has a maximum.

\begin{lemma}\label{lem1}
  The degree of $f_t$ is $mt$ if and only if $kn_{i-1}\neq e_{i-1}$,
  for $1\leq i\leq t$ and all $k\in A_i$.  Moreover, for $k\in A_i$
  such that $kn_{i-1}\neq e_{i-1}$ then
  \[\omega_{ik} = \frac{e_{i-1}}{e_{i-1}-k n_{i-1}},\]
  is a root of $f_t$.
\end{lemma}

\noindent Note that the root multiplicity of $\omega_{ij}$ corresponds to
the number of times that the value $j$ appears in $A_i$.

\begin{definition}
Let $\Omega(p)$ be the \textit{real roots of a polynomial $p$} without
repetitions, and $\Omega^+(p)$ and $\Omega^-(p)$ \textit{the set of
  positive and negative roots of $p$}, respectively. The expression
$\mu_i$ denotes \textit{the multiplicity of the root $\omega_i$} in
$\Omega(p)$.
\end{definition}

Recall from \Sref{section2} that the seed network $G_0$ is assumed to
have at least $m$ nodes, each with an in-degree of at least 1 (and
without self-loops). This assumption is key for
Lemma~\ref{lem2} to be applied. Otherwise, if there are nodes with
in-degree 0, then the minimum of the positive roots in $\Omega^+(f_t)$
would trivially be 1.

\begin{lemma}\label{lem2}
  If there is $A_i\in\A_t$ and $k\in A_i$, with
  $0<k<\floor{\frac{e_{i-1}}{n_{i-1}}}$, then $\Omega^+(f_t)$ and
  $\Omega^-(f_t)$ are non-empty sets. Moreover, $\min \Omega^+(f_t) =
  1+\frac{1}{m+\hat{m}-1}$ as $t\rightarrow\infty$.
\end{lemma}

\begin{proof}
Without loss of generality, assume that $\deg{f_t}=mt$. By Lemma
\ref{lem1}, each root of $f_t$ can be written as
\begin{equation}
\omega_{ik} = \frac{e_{i-1}}{e_{i-1}-k n_{i-1}}.\label{eq7}
\end{equation}
Since $e_{i-1}>0$, the root $\omega_k$ is positive if and only if
$e_{i-1}-k n_{i-1}> 0$. This last claim is true for some $k\in A_i$
because $0<k<\floor{\frac{e_{i-1}}{n_{i-1}}}$ by the
hypothesis. Hence, $\Omega^+(f_t)\neq \emptyset$. Moreover, if $t$ is the greatest 
index for which $k=1$ belongs to $A_t$, then the denominator of \Eref{eq7} yields the maximum positive
integer. The minimum element of $\Omega^+(f_t)$ is then given by
\begin{eqnarray*}
\min \Omega^+(f_t) &=&\frac{e_{t-1}}{e_{t-1}-n_{t-1}}\\
&=& \frac{e_{t-1}-n_{t-1}+n_{t-1}}{e_{t-1}-n_{t-1}}\\
&=& 1+\frac{n_{t-1}}{e_{t-1}-n_{t-1}}.
\end{eqnarray*}
Furthermore, note that $e_{t} = e_0+(m+\hat{m})t$ and $n_t = n_0+t$. Then
\begin{eqnarray*}
\lim_{t\rightarrow\infty}\left(\min\Omega^+(f_t)\right) &=& \lim_{t\rightarrow\infty}\left(1+\frac{n_{t-1}}{e_{t-1}-n_{t-1}}\right)\\
&=& \lim_{t\rightarrow\infty}\left(1+\frac{n_0+t-1}{e_0+(m+\hat{m})(t-1)-(n_0+t-1)}\right)\\
&=& 1+\frac{1}{m+\hat{m}-1}.
\end{eqnarray*}
\end{proof}

\begin{theorem}\label{thm1}
  If $\Omega^+(f_t)$ and $\Omega^-(f_t)$ are non-empty, and the sum of
  the multiplicities of positive roots is even, then there exist
  $a\in\real^-$, $b\in\real^+$, and $\alpha\in \real$ such
  that $\alpha\in (a, b)$ is a local maximum of the
  $G_t$-likelihood function $f_t$.
\end{theorem}

\begin{proof}
Without loss of generality, assume that there is no $k\in A_i$ such that $kn_{i-1}=e_{i-1}$ for $1\leq i \leq t$. Note that
$f_t$ can be written as the following product
\begin{eqnarray*}
f_t(x) &=& \prod\limits_{i=1}^t\prod\limits_{j=1}^{|A_i|}\frac{x \left(k n_{i-1}-e_{i-1}\right)+e_{i-1}}{e_{i-1} n_{i-1}}\\
&=& \left[\prod\limits_{i=1}^t\prod\limits_{j=1}^{|A_i|}\frac{k n_{i-1}-e_{i-1}}{e_{i-1} n_{i-1}}\right]\left[\prod\limits_{i=1}^t\prod\limits_{j=1}^{|A_i|}\left(x -\frac{e_{i-1}}{e_{i-1}-k n_{i-1}}\right)\right].
\end{eqnarray*}
Then, by letting
$C=\prod\limits_{i=1}^t\prod\limits_{j=1}^{|A_i|}\frac{k
  n_{i-1}-e_{i-1}}{e_{i-1} n_{i-1}}$ and by Lemma \ref{lem1} the
following equation holds:
\begin{equation*}
f_t(\alpha) = C \prod\limits_{i=1}^t\prod\limits_{j=1}^{|A_i|}\left(x -\omega_{ik_{ij}}\right).
\end{equation*}
The hypothesis implies that the sets of real numbers $\Omega^+(f_t)$
and $\Omega^-(f_t)$ are finite and non-empty, and hence
$\max\Omega^-(f_t)$ and $\min\Omega^+(f_t)$ exist. Let $a$
be the maximum element in $\Omega^-(f_t)$ and $b$ the
minimum element in $\Omega^+(f_t)$. Note that the function $f_t$ is
continuous and differentiable in $(a, b)$. Moreover,
since $a$ and $b$ are roots of $f_t$, $f_t(a) =
f_t(b) = 0$. By Rolle's Theorem, there exists a constant
$\alpha\in(a, b)$ such that $f^{\prime}(\alpha)=0$. That is, the
maximum argument of $f_t$ exists for $\alpha\in (a, b)$.

To show that the point $(\alpha, f_t(\alpha))$ is a local maximum of $f_t$,
consider the following two cases.

\begin{description}

\item[Case 1.] Assume that $\mu_i=1$ for each root
  $\omega_i\in\Omega(f_t)$. Since all the roots are simple, the
  function $f_t$ can be written as
\[
f_t(x) = C\prod\limits_{\omega_i\in \Omega(f_t)}\left(x-\omega_i\right).
\]
The derivative of $f_t$ with respect to $x$ is
\begin{equation}
f^{\prime}(x)=C\sum\limits_{\omega_i\in\Omega(f_t)}\prod\limits_{\omega_j\in \Omega(f_t) \atop \omega_j \neq \omega_i}\left(x-\omega_j\right). \label{eq8}
\end{equation}
Each term in the summation is a polynomial in the indeterminate
$x$ and has degree $mt-1$. Let
\begin{equation}
q_i(x)=\prod_{\omega_j\in \Omega(f_t) \atop \omega_j \neq \omega_i}\left(x-\omega_j\right).\label{eq9}
\end{equation}
By the definition of $q_i$, the root $\omega_i\in\Omega(f_t)$ does not
belong to $\Omega(q_i)$, i.e., $\Omega(q_i)\subsetneq\Omega(f_t)$. Based on \Eref{eq9}, \Eref{eq8} can be written as
\begin{equation}
f^{\prime}(x)=C \sum\limits_{\omega_i\in\Omega(f_t)}q_i(x). \label{eq10}
\end{equation}
By Lemma \ref{lem1}, there exists an in-degree $k\in A_i$ for some
$A_i\in B_t$ satisfying $\omega_{ij}=a$. In particular, there exists
$q_{s}$ such that $a$ does not belong to $\Omega(q_{s})$. Therefore,
\begin{equation}
f^{\prime}(a)=C q_s(a).\label{eq11}
\end{equation}
Since the first factor in \Eref{eq11} is positive, the sign of
$f^{\prime}(a)$ is the sign of $q_s(a)$. Note that the
product in \Eref{eq9} can be split into two factors, one
containing the positive roots of $q_s$ and another one with the
negative roots of $q_s$. More precisely,
\begin{eqnarray}
q_s(x)&=&\prod_{\omega_j\in \Omega(f_t) \atop \omega_j \neq a}\left(x-\omega_j\right)\nonumber\\
&=& \prod\limits_{\omega_j\in\Omega^+(q_s) \atop\omega_j\neq a}(x-\omega_j)\prod\limits_{\omega_j\in\Omega^-(q_s) \atop \omega_j\neq a}(x-\omega_j).\label{eq12}
\end{eqnarray}
The expression $q_s(a)$ can be rewritten as:
\begin{equation}
q_s(a)=\prod\limits_{\omega_j\in\Omega^+(q_s) \atop\omega_j\neq a}(a-\omega_j)\prod\limits_{\omega_j\in\Omega^-(q_s) \atop \omega_j\neq a}(a-\omega_j).\label{eq13}
\end{equation}
In \Eref{eq13}, by the hypothesis, the first factor is a product of an
even number of negative terms, which implies it is positive.  The
factors in the second product are positive because $a$ is the
maximum element in $\Omega^-(f_t)$ and $a < -\omega_i$, for all
$\omega_i\in\Omega^-(q_s)$.

By Lemma \ref{lem1}, there exists an in-degree $k^{\prime}\in
A_{i^\prime}$ for some $A_{i^\prime}\in B_t$ such that
$\omega_{i^{\prime}k^{\prime}}=b$; therefore, the sign of
$f^{\prime}(b)$ depends on the sign of $q_{s^{\prime}}(b)$. Using the
same argument as above:
\begin{equation}
q_{s^{\prime}}(b)=\prod\limits_{\omega_j\in\Omega^+(q_{s^{\prime}}) \atop\omega_j\neq b}(b-\omega_j)\prod\limits_{\omega_j\in\Omega^-(q_{s^{\prime}}) \atop\omega_j\neq b}(b-\omega_j).\label{eq14}
\end{equation}
Note that $b-\omega_j<0$ because $b$ is the minimum
positive root of $f_t$ and the roots are unique. Moreover, the number
of terms in the first product in \Eref{eq14} is an odd number. Thus,
this product is negative. By Lemma \ref{lem2}, all positive roots are
greater than $1$, so the second product in \Eref{eq14} is
positive. Consequently, $q_{s^{\prime}}(b)$ is
negative. Therefore, $f^{\prime}(a)>0$ and
$f^{\prime}(b)<0$. Hence, $(\alpha,f(\alpha))$ is a local maximum of
$f_t$.

\item[Case 2.]

  Assume that there is at least one root $\omega_i$ in $\Omega(f_t)$
  with $\mu_i \geq 2$. Since $\mu_i\geq 1$ for all roots
  $\omega_i\in\Omega(f_t)$, the function $f_t$ can be written as
  \[
f_t(x) = C\prod\limits_{\omega_i\in \Omega(f_t)}\left(x-\omega_{i}\right)^{\mu_i}.
\]
The derivative of $f_t$ with respect to $x$ is
\begin{equation}
f^{\prime}(x)= C \sum\limits_{\omega_i\in \Omega(f_t)}\left[\mu_i(x-\omega_{i})^{\mu_i-1}\prod\limits_{\omega_j\in \Omega(f_t)\atop \omega_j\neq \omega_i}\left(x -\omega_i\right)^{\mu_i}\right]. \label{eq15}
\end{equation}
The summands in \Eref{eq15} are polynomials, each with the same roots
as $f_t$, but with multiplicity of at most $\mu_i$. As for case 1, the
goal is to show that $f^{\prime}(a+\epsilon)>0$ and
$f^{\prime}(b-\epsilon)<0$, for some $\epsilon>0$.\\
By defintion the elements of $\Omega(f_t)$ are unique, then it can be indexed as
$\Omega(f_t)=~\{\omega_{i_1},\omega_{i_2},\ldots,\omega_{i_s}\}$. Expanding the second factor in
\Eref{eq15}:
\begin{eqnarray}
\sum\limits_{\omega_i\in \Omega(f_t)}&&\left[\mu_i(x-\omega_{i})^{\mu_i-1}\prod\limits_{\omega_j\in \Omega(f_t)\atop \omega_j\neq \omega_i}\left(x -\omega_i\right)^{\mu_i}\right] \\
&=& \mu_{i_1}(x-\omega_{i_1})^{\mu_{i_1}-1}[\left(x -\omega_{i_2}\right)^{\mu_{i_2}} \left(x -\omega_{i_3}\right)^{\mu_{i_3}} \cdots ]\nonumber\\
&&+\mu_{i_2}(x-\omega_{i_2})^{\mu_{i_2}-1}[\left(x -\omega_{i_1}\right)^{\mu_{i_1}} \left(x -\omega_{i_3}\right)^{\mu_{i_3}} \cdots ]\nonumber\\
&&+\mu_{i_3}(x-\omega_{i_3})^{\mu_{i_3}-1}[\left(x -\omega_{i_1}\right)^{\mu_{i_1}} \left(x -\omega_{i_2}\right)^{\mu_{i_2}} \cdots ] \nonumber\\
&&+\cdots \nonumber\\
&=& \left[\prod\limits_{\omega_i\in \Omega(f_t)}(x-\omega_i)^{\mu_i-1}\right]\sum\limits_{\omega_i\in \Omega(f_t)}\left[\mu_i\prod\limits_{\omega_j\in \Omega(f_t)\atop \omega_j\neq \omega_i}(x-\omega_j)\right].
\label{eq16}
\end{eqnarray}
Using \Eref{eq16}, \Eref{eq15} can be written as 
\begin{equation}
f^{\prime}(x)=C\left[\prod\limits_{\omega_i\in \Omega(f_t)}(x-\omega_i)^{\mu_i-1}\right]\sum\limits_{\omega_i\in \Omega(f_t)}\left[\mu_i\prod\limits_{\omega_j\in \Omega(f_t)\atop \omega_j\neq \omega_i}(x-\omega_j)\right].\label{eq17}
\end{equation}

In \Eref{eq17}, the factor $\sum\limits_{\omega_i\in \Omega(f_t)}\left[\mu_k\prod\limits_{\omega_j\in \Omega(f_t)\atop \omega_j\neq \omega_i}(x-\omega_j)\right]$ is equal to
$\sum\limits_{\omega_i\in \Omega(f_t)}\mu_i q_i(x)$. On the one hand,
$q_i(a) > 0$ if and only if $|\Omega(f_t)|$ is even. Since
$q_i$ is continuous over its domain, there exist $M>0$ and $\delta>0$
such that if
$(a+\epsilon)\in\:(a-\epsilon,b+\epsilon)$ and
$|(a+\epsilon)-a|<\delta$, then
$q_i(a+\epsilon)>0$~\cite{Apostol1974}. The second term in
\Eref{eq17} can be rewritten as
\begin{equation}
\prod\limits_{\omega_i\in\Omega^-(f_t)}(x-\omega_i)^{\mu_i-1}\prod\limits_{\omega_i\in\Omega^+(f_t)}(x-\omega_i)^{\mu_i-1}.\label{eq18}
\end{equation}
The first term in \Eref{eq18} is positive evaluated at
$a+\epsilon$. The second term has the property that the sum of
all exponents is an even number and, evaluated at $a+\epsilon$,
is also positive. Hence, $f^\prime(a+\epsilon)>0$.

If $|\Omega(f_t)|$ is odd, the function $q_i$ evaluated at $a$
is negative. Using the same argument as above, it can be shown that
$q_i(a+\epsilon)<0$. There exists an odd number of even
multiplicities by the hypothesis. Consequently, the second factor in
\Eref{eq18}, evaluated at $a+\epsilon$, is negative. In either
case $f^\prime(a+\epsilon)>0$. Similarly, it can be shown that
$f^\prime(b-\epsilon)<0$.
\end{description}

\noindent Therefore, $(\alpha,f(\alpha))$ is a local
maximum of $f_t$.
\end{proof}

\begin{figure}[htbp]
\centering
\includegraphics[width=200pt]{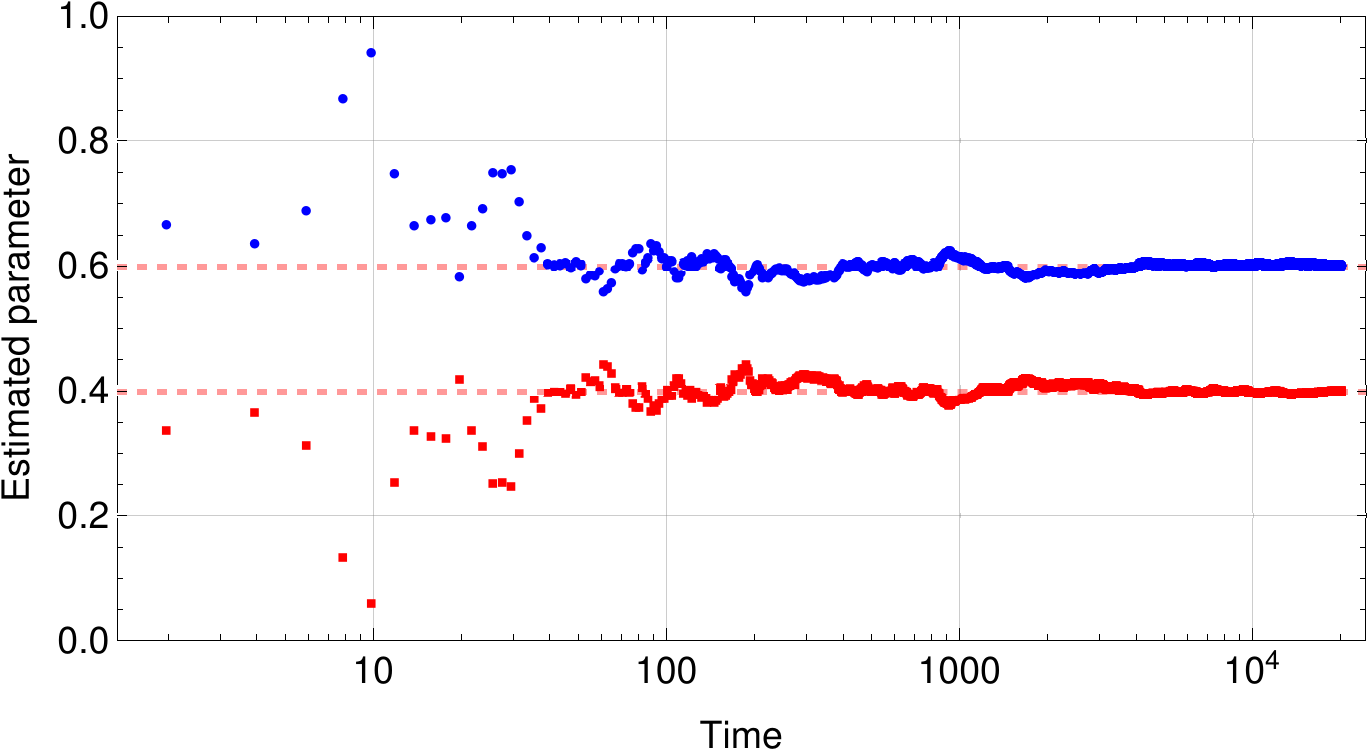}
\caption{Convergence of the parameters associated to preferential
  (dots) and random (filled squares) attachment mechanisms by
  application of Theorem~\ref{thm1} for the complete graph $G_0$ with
  $N_0=3$ nodes, and a sequence $(G_t)_{t\leq 20000}$ generated by the
  algorithm in Definition~\ref{def.netmodel.netmodel} with parameters
  $m = 5$, $\hat{m} = 3$, and $\alpha=0.6$.}
\label{fig:test}
\end{figure}

To illustrate the application of Theorem~\ref{thm1}, consider the
complete graph $G_0$ with $N_0=3$ nodes, and a sequence $(G_t)_{t\leq
  20000}$ generated by the algorithm in
Definition~\ref{def.netmodel.netmodel} with parameters $m = 5$,
$\hat{m} = 3$, and $\alpha=0.6$. At each time $t$, the
$G_t$-likelihood function is the one in
Definition~\ref{def.maxlike.likelihood}. An estimate for $\alpha$ is
computed maximizing the $G_t$-likelihood function under the conditions
of Theorem~\ref{thm1}.  \Fref{fig:test} shows how maximizing the
$G_t$-likelihood function yields the estimate $\tilde{\alpha}$ for the
parameter~$\alpha$.

\section{Using an Expectation-Maximization Algorithm}
\label{sec.em}

This section presents how determining the contribution of the two
attachment mechanisms used by the algorithm in
Definition~\ref{def.netmodel.netmodel} is equivalent to computing the
likelihood estimator with an Expectation-Maximization (EM) algorithm.
The basic idea is to associate a complete-data problem, which is
better suited for maximum likelihood estimation, to a given
incomplete-data problem, for which the same estimation can become a
wild-goose chase. For an EM algorithm to be effective, two important
requirements need to be met: first, it needs to be proved convergent;
second, it needs to be efficient. This section addresses these two
requirements by identifying sufficient conditions for the algorithm to
converge and by presenting a recursive definition for estimating
$\alpha$, which can be used for incremental computation in the EM
algorithm.

\begin{definition}\label{def.em.pmf}
  Let $X_t$ be a random variable that characterizes the in-degree of
  the selected nodes due to an attachment mechanism up to time $t$
  from $G_t$. Moreover, let $\mathbb{P}(x\msp\theta)$ denote the
  probability mass function (pmf) as a function of the parameter
  vector $\theta$.
\end{definition}
Recall that, for a fixed time $t$, \Eref{eq6} defines the likelihood
of forming an edge $(u,v)$ from the new node $u$ to node $v\in
V_{t-1}$ through an attachment mechanism. Such an equation can be
rewritten as:
\begin{eqnarray}
f_t(\alpha)&= \prod\limits_{i=1}^t\prod\limits_{j=1}^{|A_i|}\left[\alpha\frac{k_{ij}}{e_{i-1}}+(1-\alpha)\frac{1}{n_{i-1}}\right].
\label{eq19}
\end{eqnarray}
The set $\A_t=\{A_1,\ldots,A_t\}$ is a collection of samples of length
$m$ generated by the random variable $X_t$ (see
Definition~\ref{def.em.pmf}).

The EM algorithm computes iteratively a
maximum likelihood estimator for data with unobserved
variables~\cite{Dempster1977,bordes2007algorithm}. In this case, the input to the EM
algorithm comes from a mixed distribution in which the mixture weight
is unknown. Based on \Eref{eq19}, the incomplete likelihood function
can be modeled as follows:
  \begin{eqnarray}
\mathcal{L}(\alpha\msp\A_t)= \sum\limits_{i=1}^{t}\sum\limits_{j=1}^{|A_i|}\log\left(\alpha\frac{k_{ij}}{e_{i-1}}+(1-\alpha)\frac{1}{n_{i-1}}\right).
 \label{eq20}
  \end{eqnarray}
where $k_{ij}$ is the $j$-th in-degree in the sample $A_i$ (which, as
assumed in Section~\ref{sec.maxlike} and without loss of generality,
is sorted in ascending order) and $|A_i|$ is the length of $A_i$.  The
attachment mechanism operates on two groups of nodes. One group
consists of nodes that attach to the network using preferential
attachment and the other of the nodes that attach using random
attachment. Consider a binary variable $Y$ for each occurrence of an
in-degree $k$ to indicate whether the observation has been selected by
preferential attachment or by random attachment. For $k_{ij}$, define
$y_{ij}=1$ if $k_{ij}$ is in group one and $y_{ij}=0$
otherwise~\cite{Mclachlan2007}.

Notice that $\mathbb{P}(y_{ij}=1\msp\alpha)=\alpha$ and
$\mathbb{P}(y_{ij}=0\msp\alpha)=1-\alpha$. The conditional probabilities of the data
and the unobservable data are:
\begin{eqnarray*}
\mathbb{P}(k_{ij}\msp y_{ij},\alpha)&=&\left(\frac{k_{ij}}{e_{i-1}}\right)^{y_{ij}}\left(\frac{1}{n_{i-1}}\right)^{1-y_{ij}}\\
\mathbb{P}(y_{ij}\msp \alpha) &=& \alpha^{y_{ij}}(1-\alpha)^{1-y_{ij}}.
\end{eqnarray*}
By defining $Z=(\A_t,Y)$ as the complete-data, the complete
log-likelihood function is given by
\begin{eqnarray}
\mathcal{L}_c(Z\msp \alpha)&=& \sum\limits_{i=1}^{t}\sum\limits_{j=1}^{|A_i|}\log\left(\mathbb{P}(k_{ij}\msp y_{ij},\alpha)\mathbb{P}(y_{ij}\msp \alpha\right)\nonumber\\
&=& \sum\limits_{i=1}^{t}\sum\limits_{j=1}^{|A_i|} \log \left(\left[\alpha\frac{k_{ij}}{e_{i-1}}\right]^{y_{ij}}\left[(1-\alpha)\frac{1}{n_{i-1}}\right]^{1-y_{ij}}\right)\nonumber\\
&=& \sum\limits_{i=1}^{t}\sum\limits_{j=1}^{|A_i|}y_{ij}\log \left(\alpha\frac{k_{ij}}{e_{i-1}}\right)\nonumber\\
&&+\sum\limits_{i=1}^{t}\sum\limits_{j=1}^{|A_i|}\left(1-y_{ij}\right)\log\left((1-\alpha)\frac{1}{n_{i-1}}\right).
\label{eq21}
\end{eqnarray}
\Eref{eq21} is linear in the unobservable data $y_{ij}$. Consider a
function $Q(\alpha\msp\alpha^{(d)})$ that represents the conditional
expectation given the observed data using the $d$-th fit for the
unknown parameter $\alpha$. In particular,
\begin{eqnarray}
Q\left(\alpha\msp\alpha^{(d)}\right) &=& \sum\limits_{i=1}^{t}\sum\limits_{j=1}^{|A_i|}\mathbb{P}\left(y_{ij}=1\msp k_{ij},\alpha^{(d)}\right)\log \left(\alpha\frac{k_{ij}}{e_{i-1}}\right)\nonumber\\
&&+\sum\limits_{i=1}^{t}\sum\limits_{j=1}^{|A_i|}\mathbb{P}\left(y_{ij}=0\msp k_{ij},\alpha^{(d)}\right)\log\left((1-\alpha)\frac{1}{n_{i-1}}\right).
\label{eq22}
\end{eqnarray}

By Bayes' Theorem,
\begin{eqnarray}
\mathbb{P}\left(y_{ij}=1\msp k_{ij},\alpha^{(d)}\right) &=& \frac{\mathbb{P}\left(k_{ij}\msp y_{ij}=1,\alpha^{(d)}\right)\mathbb{P}\left(y_{ij}=1\msp \alpha^{(d)}\right)}{\mathbb{P}\left(k_i\msp \alpha^{(d)}\right)}\nonumber\\
&=&  \frac{\frac{k_{ij}}{e_{i-1}}\alpha^{(d)}}{\frac{k_{ij}}{e_{i-1}}\alpha^{(d)}+\frac{1}{n_{i-1}}(1-\alpha^{(d)})}
\label{eq23}
\end{eqnarray}
And
\begin{eqnarray}
\mathbb{P}\left(y_{ij}=0\msp k_{ij},\alpha^{(d)}\right) &=& 1-\mathbb{P}\left(y_{ij}=1\msp k_{ij},\alpha^{(d)}\right)\nonumber\\
&=&  \frac{\frac{1}{e_{i-1}}\left(1-\alpha^{(d)}\right)}{\frac{k_{ij}}{e_{i-1}}\alpha^{(d)}+\frac{1}{n_{i-1}}(1-\alpha^{(d)})}
\label{eq24}
\end{eqnarray}
Consider $x_1=\alpha$ and $x_2=1-\alpha$. In order to maximize
$Q\left(\alpha\msp\alpha^{(d)}\right)$, the Lagrange multiplier $\lambda$ 
with the constrain $x_1 + x_2 = 1$ is introduced with the goal of solving the
following equation:
\begin{equation}
\sum\limits_{i=1}^{t}\sum\limits_{j=1}^{|A_i|} \frac{\mathbb{P}\left(k_{ij}\msp y_{ij}=1,\alpha^{(d)}\right)}{x_1}-\lambda=0
\label{eq25}
\end{equation} 
\Eref{eq25} has the solution:
\begin{equation}
\alpha^{(d+1)} = \frac{1}{\sum\limits_{i=1}^{t}|A_i|}\sum\limits_{i=1}^{t}\sum\limits_{j=1}^{|A_i|}\frac{\frac{k_{ij}}{e_{i-1}}\alpha^{(d)}}{\frac{k_{ij}}{e_{i-1}}\alpha^{(d)}+\frac{1}{n_{i-1}}(1-\alpha^{(d)})}.
\label{eq26}
\end{equation}

\Eref{eq26} provides a recursive formulation of $\alpha$ that can be
used for estimation purposes in Algorithm~\ref{alg:em}, namely, in the
EM algorithm. This algorithm computes an estimate for $\alpha$ on an
input column block matrix $[A_1\:A_2\:\cdots A_t]^T$ where $A_i\in\mathbb{Z}^{m\times 3}$. For a fixed
$A_i$, the values $a_{i1}$, $a_{i2}$, and $a_{i3}$ denote the
in-degree of the target node, the number of edges, and the number of
nodes in the network, respectively. Line 3 computes the probability
of the new edge forming by preferential attachment or random
attachment mechanism (it implements the E-step of the EM
algorithm). Lines 2-6 approximate the value for the parameter
$\alpha$ until a sufficiently accurate value is reached. Note that the recursive
definition of $\alpha$ presented in \Eref{eq26} is used in line 4 (it
implements the M-step of the EM algorithm). Finally, 
the convergence of Algorithm~\ref{alg:em} is obtained as a
corollary of Lemma 2.1 in~\cite{diebolt-em-93}.

\begin{algorithm}[htbp]
\caption{EM Algorithm}\label{alg:em}
\begin{algorithmic}[1]
\renewcommand{\algorithmicrequire}{\textbf{Input:}}
\renewcommand{\algorithmicensure}{\textbf{Output:}}
\Require {A $t$-column block matrix $A$ of integers with $t > 0$ and an error bound $\epsilon > 0.$}
\Ensure An estimate for $\alpha$.\\
$\alpha_c$ = $0.5$;
$\alpha_b$ = $0$;
\While {$|\alpha_c-\alpha_b|\geq \epsilon$}\\
\quad $\alpha_b\leftarrow\alpha_c$\\
\quad $e\leftarrow \sum\limits_{i=1}^{t}\sum\limits_{j=1}^{|A_i|}\frac{\frac{A_{ij1}}{A_{ij2}}\alpha_c}{\frac{A_{ij1}}{A_{ij2}}\alpha_c+\frac{1}{A_{ij3}}(1-\alpha_c)}$\\
\quad $\alpha_c \leftarrow \frac{e}{\sum_{i=1}^{t}|A_i|}$

\EndWhile\\
\Return $\alpha_c$
\end{algorithmic}
\end{algorithm}

\begin{corollary}\label{cor.em.emconv}
  Algorithm~\ref{alg:em} converges.
\end{corollary}

The main observation in the proof of Corollary~\ref{cor.em.emconv} is
that the function $\log \mathbb{P}(\alpha\msp\A_t)$ in \Eref{eq20} is
concave.

\begin{figure}[htbp]
\centering
\includegraphics[width=200pt]{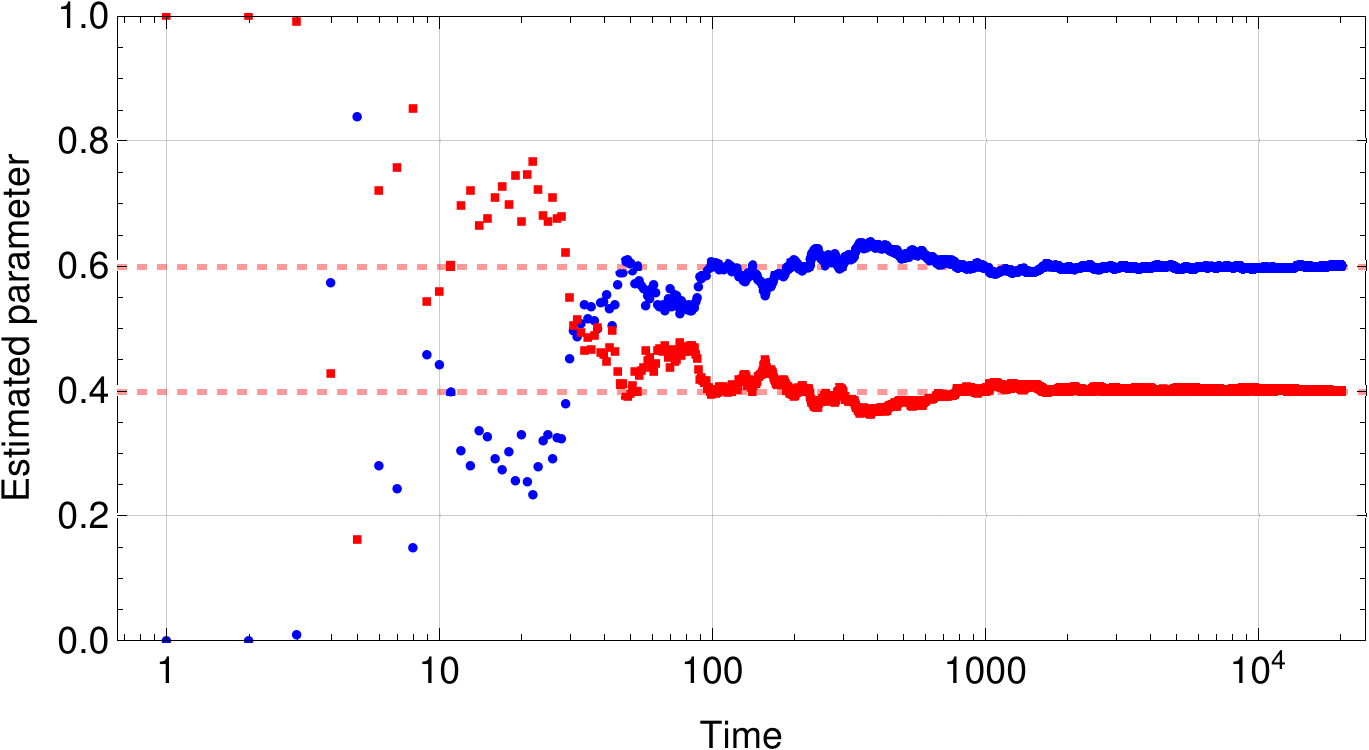}
\caption{Convergence of the parameters associated to preferential
  (dots) and random (filled squares) attachment mechanisms by the EM
  Algorithm for the complete graph $G_0$ with $N_0=3$ nodes, and a
  sequence $(G_t)_{t\leq 20000}$ generated by the algorithm in
  Definition~\ref{def.netmodel.netmodel} with parameters $m = 5$,
  $\hat{m} = 3$, and $\alpha=0.6$.}
\label{fig:testem}
\end{figure}

Recall the complete graph $G_0$ with $N_0=3$ nodes, and a sequence
$(G_t)_{t\leq 20000}$ generated by the algorithm in
Definition~\ref{def.netmodel.netmodel} with parameters $m = 5$,
$\hat{m} = 3$, and $\alpha=0.6$. An estimate for $\alpha$ computed by
Algorithm~\ref{alg:em} (i.e., the EM Algorithm) is depicted in
Figure~\ref{fig:testem}. The final estimate computed for $\alpha$ by
the EM Algorithm closely approximates the one computed by using
Theorem~\ref{thm1} at the end of Section~\ref{sec.maxlike}.

\section{In-degree Distribution}
\label{section4}

This section characterizes the in-degree distribution of the nodes in
a sequence of networks $(G_t)_{t \leq T}$ generated by the algorithm
presented in Definition~\ref{def.netmodel.netmodel}, for some $T \in
\nats$, from a given seed network $G_0$. It also presents how the
dynamics of the in-degree distribution converges to a stationary
distribution and illustrates the approach with experiments on
sequences of networks $(G_t)_{t \leq T}$.

\begin{definition}\label{def.indegree.pmf}
  Let $K_t$ be a random variable that characterizes the in-degree of a
  node selected uniformly at random at time $t$ from $G_t$.  Moreover,
  let $P_t(k)=\mathbb{P}(K_t=k)$ denote the probability that $K_t$ is
  equal to $k$ at time $t$.
\end{definition}

The notion of asymptotic equivalence between two real sequences is
used to prove the existence of
$\lim_{t\rightarrow\infty}P_t(k)$. Theorem \ref{thm2} ensures that the
probability of the in-degree distribution of the model converges.

\begin{theorem}\label{thm2}
  For $t \in \nats$, $\lim_{t\rightarrow\infty}P_t(k)$ exists for all
  $k\geq \hat{m}$.
\end{theorem}

\begin{proof}
  First, note that the expected number of nodes of in-degree $k \geq \hat{m}$ is
  \begin{eqnarray}
    n_{t}P_{t}(k) &=& n_{t-1}P_{t-1}(k)-m\pi_t(k\mid\alpha)n_{t-1}P_{t-1}(k)\nonumber\\
    &&+m\pi_t(k-1\mid\alpha)n_{t-1}P_{t-1}(k-1). \label{eq27}
  \end{eqnarray}
  
According to \Eref{eq5}
\begin{eqnarray}
n_{t}P_{t}(k)&=& \left(n_{t-1}-\frac{m\alpha k n_{t-1}}{e_{t-1}}-m(1-\alpha)\right)P_{t-1}(k)\nonumber\\
&&+\left(\frac{m\alpha (k-1)n_{t-1}}{e_{t-1}}+m(1-\alpha)\right)P_{t-1}(k-1). \label{eq28}
\end{eqnarray}

Since $\hat{m}$ nodes establish an edge to the new node, the expected
number of nodes of in-degree $k=\hat{m}$ is
\begin{equation}
n_{t}P_{t}(\hat{m})= \left(n_{t-1}-\frac{m\hat{m}\alpha n_{t-1}}{e_{t-1}}-m(1-\alpha)\right)P_{t-1}(\hat{m})+1.\label{eq29}
\end{equation}   
The proof that the limit exists follows by induction on $k$.

\begin{description}

\item[Base case.] When $k=\hat{m}$, by using \Eref{eq29},
  $P_t(\hat{m})$ can be expressed using the recurrence
  \begin{equation}
P_t(\hat{m})=\frac{1}{n_t}\left(n_{t-1}-\frac{m\hat{m}\alpha n_{t-1}}{e_{t-1}}-m(1-\alpha)\right)P_{t-1}(\hat{m})+\frac{1}{n_t}.\label{eq30}
  \end{equation}
  This is a non-autonomous, first-order difference equation. It can be
  shown by induction over $t$ that the solution of \Eref{eq30} is
  given by
  \begin{equation}
P_t(\hat{m})=\prod\limits_{i=1}^{t}a_iP_0(\hat{m})+\sum\limits_{i=1}^{t}\left[\prod\limits_{j=i+1}^{t}a_i\right]b_i,\label{eq31}
  \end{equation}
  where $a_t=\frac{1}{n_t}\left(n_{t-1}-\frac{m\hat{m}\alpha
    n_{t-1}}{e_{t-1}}-m(1-\alpha)\right)$ and $b_t=
  \frac{1}{n_t}$. The first term in \Eref{eq31}, can be written in
  terms of the Gamma functions as
  \[\prod\limits_{i=1}^{t}a_iP_0(\hat{m}) = \frac{\Gamma(i-\xi_1)\Gamma(i-\xi_2)P_0(\hat{m})}{\Gamma(\frac{e_0}{m+\hat{m}}+i)\Gamma(n_0+i)},\]
  where $\xi_1$ and $\xi_2$ are constant real numbers that do not
  depend on time. Moreover,
\[\lim\limits_{t\rightarrow\infty}\frac{\Gamma(t-\xi_1)\Gamma(t-\xi_2)P_0(\hat{m})}{\Gamma(\frac{e_0}{m+\hat{m}}+t)\Gamma(n_0+t)}=0.\]

It can further be shown that the second term in \Eref{eq31} is the
convergent series
\[\sum\limits_{i=1}^{\infty}\left[\prod\limits_{j=i+1}^{\infty}a_i\right]b_i=\frac{m + \hat{m}}{m + \hat{m} + m^2 + m \hat{m} - m^2 \alpha}.\]

Therefore,
\[\lim\limits_{t\rightarrow\infty}P_t(\hat{m})=\frac{m + \hat{m}}{m + \hat{m} + m^2 + m \hat{m} - m^2 \alpha}.\]

\item[Inductive step.]  Let $k>\hat{m}$ and assume that
  $\lim_{t\rightarrow\infty}P_t(k)$ exists for all $k>\hat{m}$. For a
  large enough $t$, $P_{t-1}(k)\sim P_t(k)$ and $P_{t-1}(k+1)\sim
  P_t(k+1)$. By \Eref{eq30} and following a line of argument similar
  to the one in~\cite[Theorem~1]{Ruiz2017}:
\begin{eqnarray*}
\left(1+\frac{m\alpha (k+1)}{m+\hat{m}}+m(1-\alpha)\right)P_t(k+1) \sim \left(\frac{m\alpha k}{m+\hat{m}}+m(1-\alpha)\right)P_t(k).
\end{eqnarray*}

By inductive hypothesis
\[\lim_{t\rightarrow\infty}P_t(k+1)=\frac{\frac{m\alpha k}{m+\hat{m}}+m(1-\alpha)}{1+\frac{m\alpha (k+1)}{m+\hat{m}}+m(1-\alpha)}\lim_{t\rightarrow\infty}P_t(k).\]
\end{description}
Therefore, $\lim_{t\rightarrow\infty}P_t(k)$ exists for all $k\geq \hat{m}$.
\end{proof}

\Eref{eq28} indicates that the expected number of nodes of in-degree
$k\geq \hat{m}$ is equal to the difference between the expected number
of nodes of in-degree $k$ selected at time $t-1$ by the attachment
process and the expected number of nodes of in-degree $k-1$ that
establish an edge with the new node.

Corollary~\ref{cor.indegree.behavior} characterizes the in-degree
distribution of the network.

\begin{corollary}\label{cor.indegree.behavior}
If $k \geq \hat{m}$, then the asymptotic behavior of the expected
complementary cumulative in-degree distribution satisfies
\[
\bar{F}_\infty(k)= \left\{ \begin{array}{lll}
             \left(\frac{m}{1+m}\right)^{k-\hat{m}} &   ,  & \alpha = 0 \\
             \\ \frac{\Gamma\left(\hat{m}+m\hat{m}+m\right)\Gamma(k)}{\Gamma(\hat{m})\Gamma\left(k+\frac{m+\hat{m}}{m}\right)} &  , & \alpha=1 \wedge\hat{m}\geq 1 \\
             \\ \frac{\Gamma \left(\frac{\hat{m}+m (1+m+\hat{m}-\alpha  m)}{m \alpha }\right) \Gamma \left(k+\frac{(m+\hat{m}) (1-\alpha)}{\alpha }\right)}{\Gamma \left(\frac{\hat{m}+m(1-\alpha)}{\alpha }\right) \Gamma \left(\frac{\hat{m}+m (m+\hat{m}+k \alpha -(m+\hat{m}) \alpha +1)}{m \alpha }\right)} &  ,  & 0<\alpha<1
             \end{array}
   \right.
   \]
\end{corollary}

\begin{proof}
Let $P_\infty(k)$ denote the limit of $P_t(k)$ as $t$ tends to
infinity. According to Theorem~\ref{thm2}, $P_\infty(k)$ can be
written as
\begin{equation}
P_\infty(k)= \left\{ \begin{array}{lll}
             \frac{m+\hat{m}}{m^2+m \hat{m}+m+\hat{m}-\alpha  m^2} &   ,  & k = \hat{m} \\
             \\ \frac{\alpha  \left(k m-m^2-m \hat{m}-m\right)+m^2+m \hat{m}}{\alpha  \left(k m-m^2-m \hat{m}\right)+m^2+m \hat{m}+m+\hat{m}}P_{\infty}(k-1)&  ,  & k> \hat{m}
             \end{array}
   \right.\label{eq32}
\end{equation}
\Eref{eq32} defines a recurrence relation that varies as a function of the value of
$\alpha$. For $\alpha = 0$, the solution of the recurrence is
\begin{equation}
P_\infty(k) = \frac{1}{m+1}\left(\frac{m}{m+1}\right)^{k-\hat{m}}.\label{eq33}
\end{equation}  
For $\alpha=1$ and $\hat{m}\geq 1$, the solution is given by
\begin{equation}
P_\infty(k) = \frac{(m+\hat{m})\Gamma(\hat{m}+m\hat{m}+m)\Gamma(k)}{m\Gamma(\hat{m})\Gamma\left(k+\frac{\hat{m}+2m}{m}\right)}.\label{eq34}
\end{equation}
Furthermore, for $0<\alpha<1$, the recurrence has solution
\begin{equation}
P_\infty(k) =\frac{(m+\hat{m}) \Gamma \left(\frac{r+m (1+m+\hat{m}- m\alpha)}{m\alpha }\right) \Gamma \left(k+\frac{(m+\hat{m}) (1-\alpha)}{\alpha }\right)}{m \alpha\: \Gamma \left(\frac{m+\hat{m}-m\alpha  }{\alpha }\right) \Gamma \left(\frac{\hat{m}+m (m+\hat{m}+k \alpha -(m+\hat{m}) \alpha +\alpha +1)}{m \alpha }\right)}.\label{eq35}
\end{equation}
Since $\bar{F}_\infty(k)=\mathbb{P}[K\geq
  k]=1-\sum_{j=\hat{m}}^{k-1}P_\infty(j)$, by using \Eref{eq33},
\Eref{eq34} and \Eref{eq35}, the desired result is obtained.
\end{proof}

\begin{figure}[hbt!]
    \centering
    \subfloat[][]{\includegraphics[width=170pt]{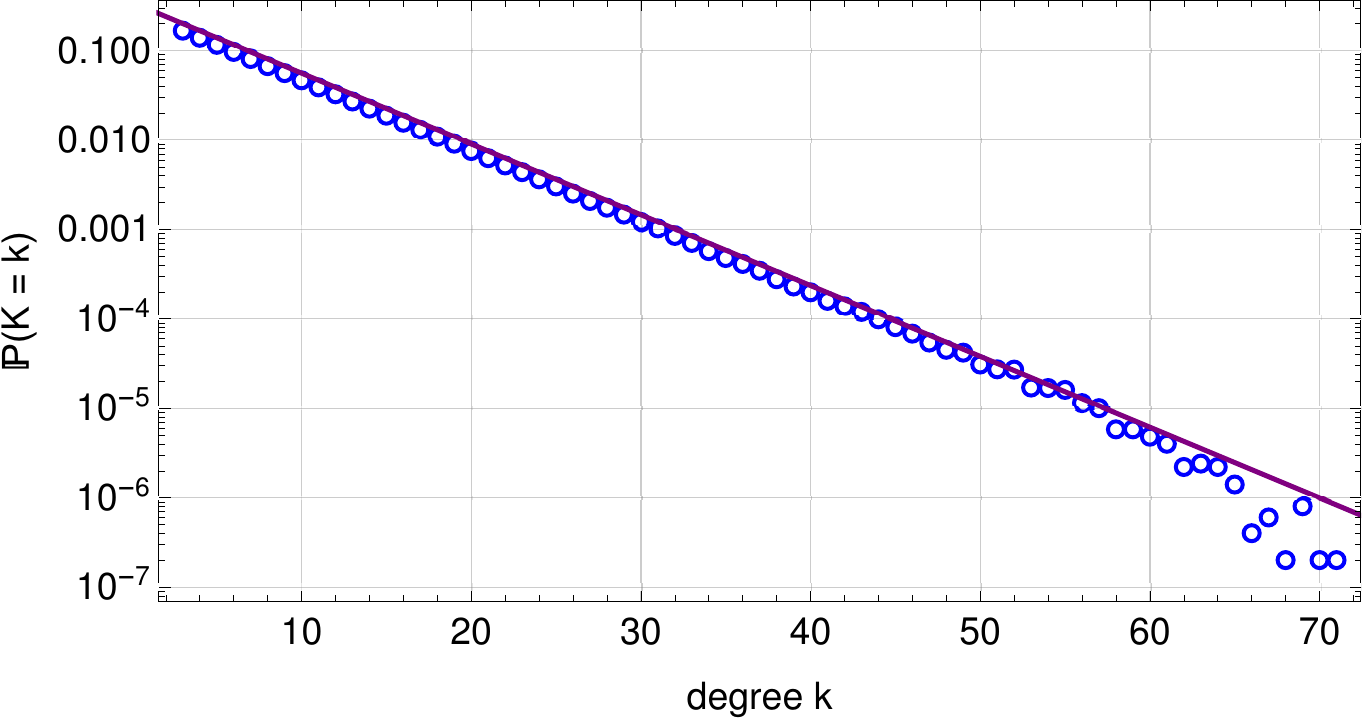}}~
    \subfloat[][]{\includegraphics[width=170pt]{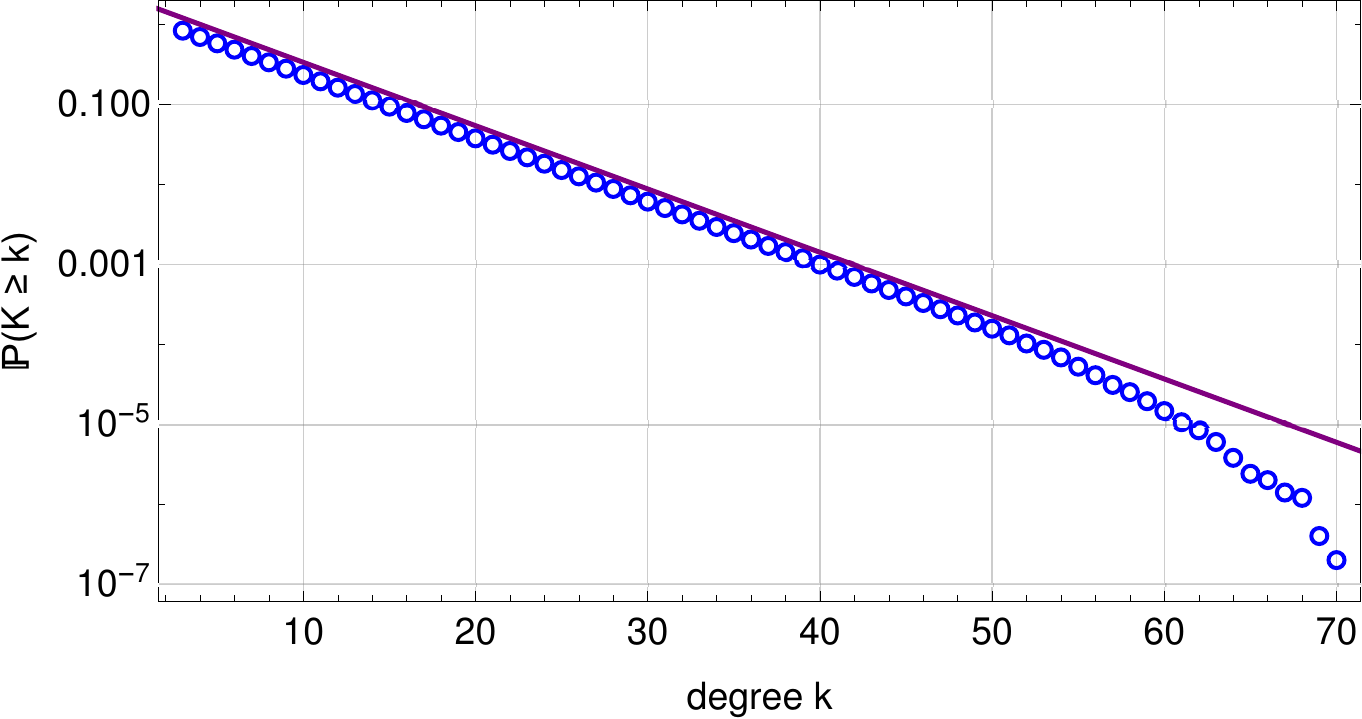}}

    \subfloat[][]{\includegraphics[width=170pt]{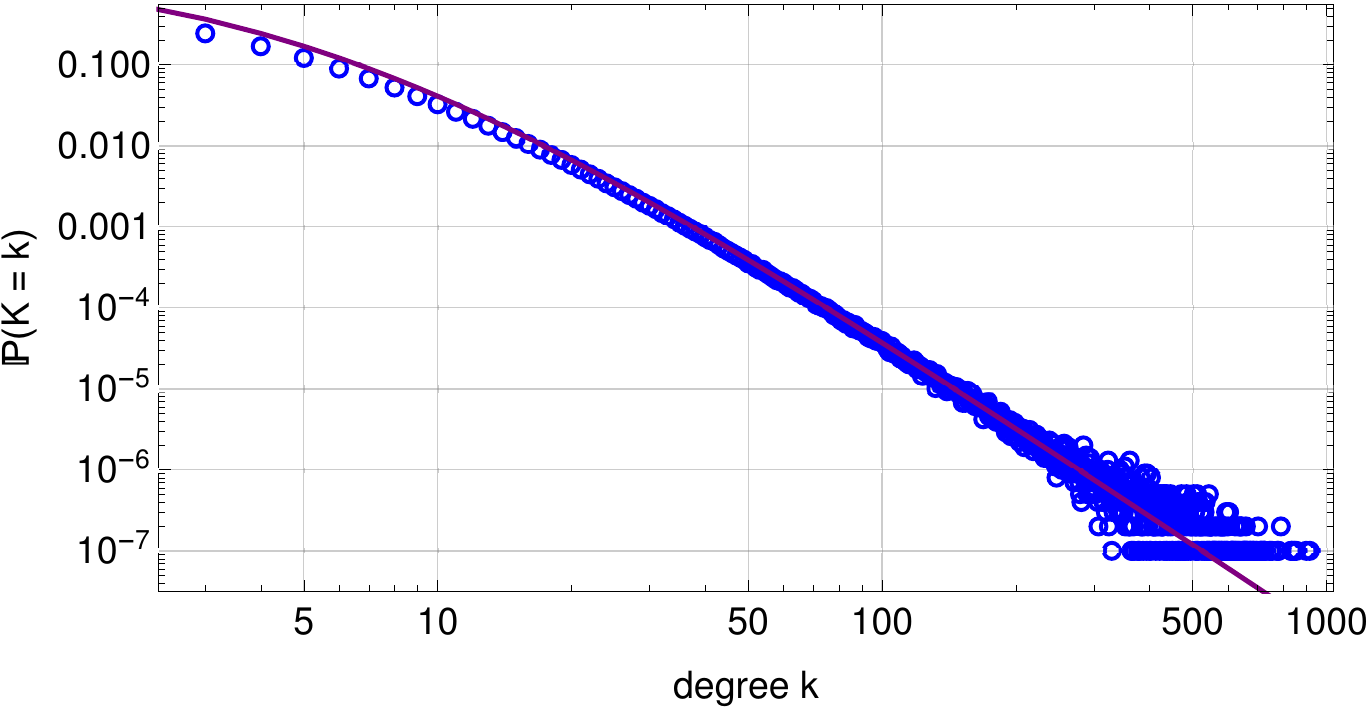}}~
    \subfloat[][]{\includegraphics[width=170pt]{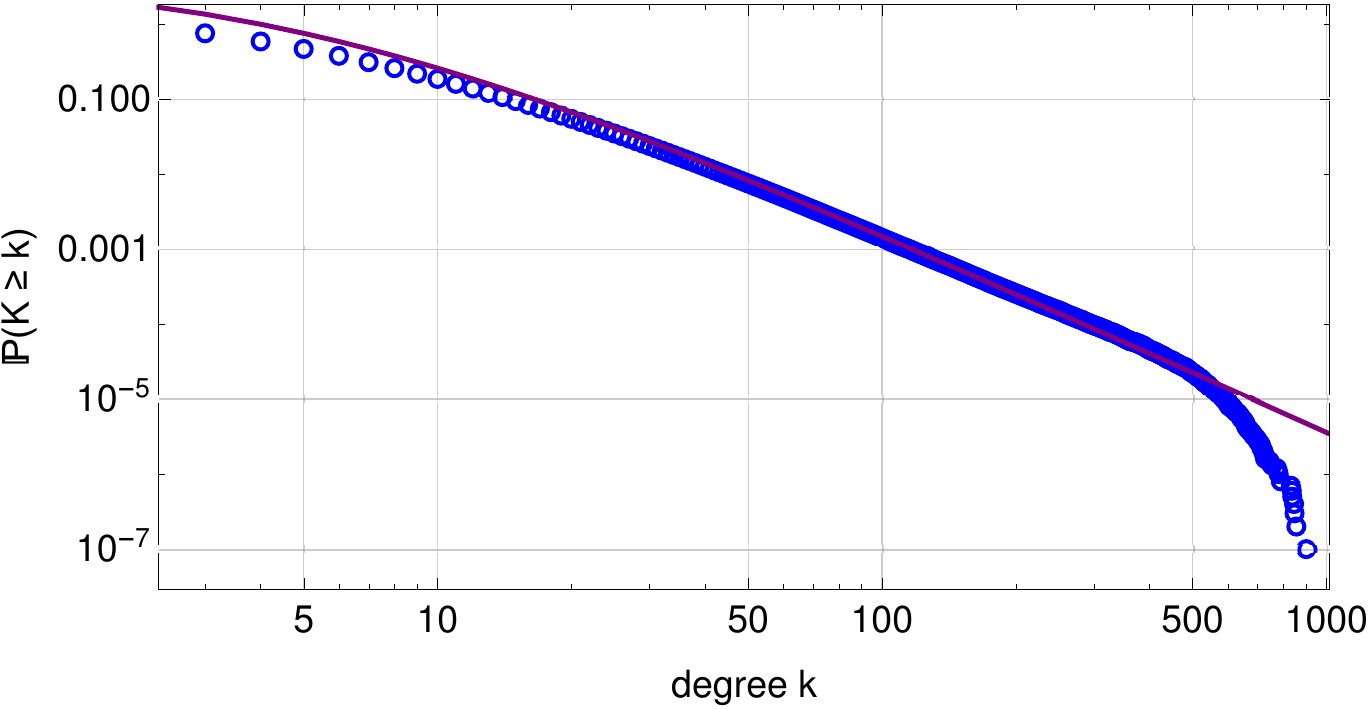}}

    \subfloat[][]{\includegraphics[width=170pt]{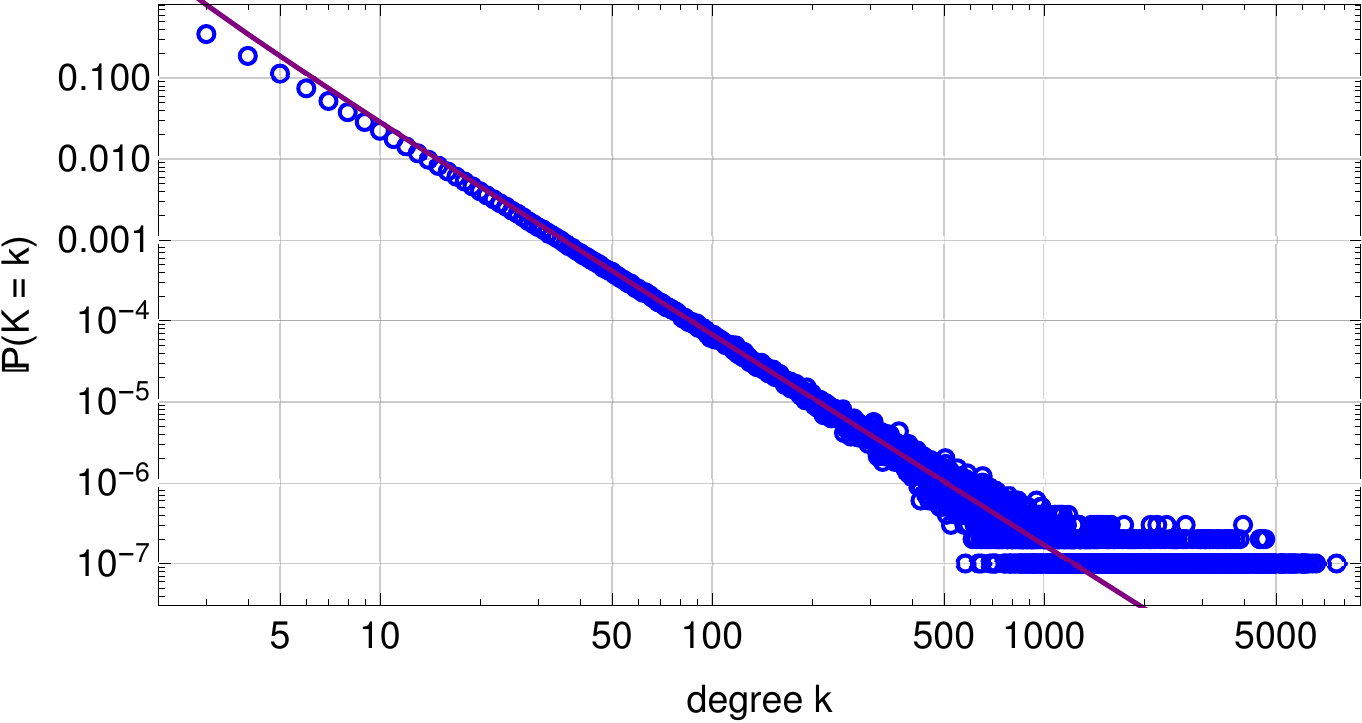}}~
    \subfloat[][]{\includegraphics[width=170pt]{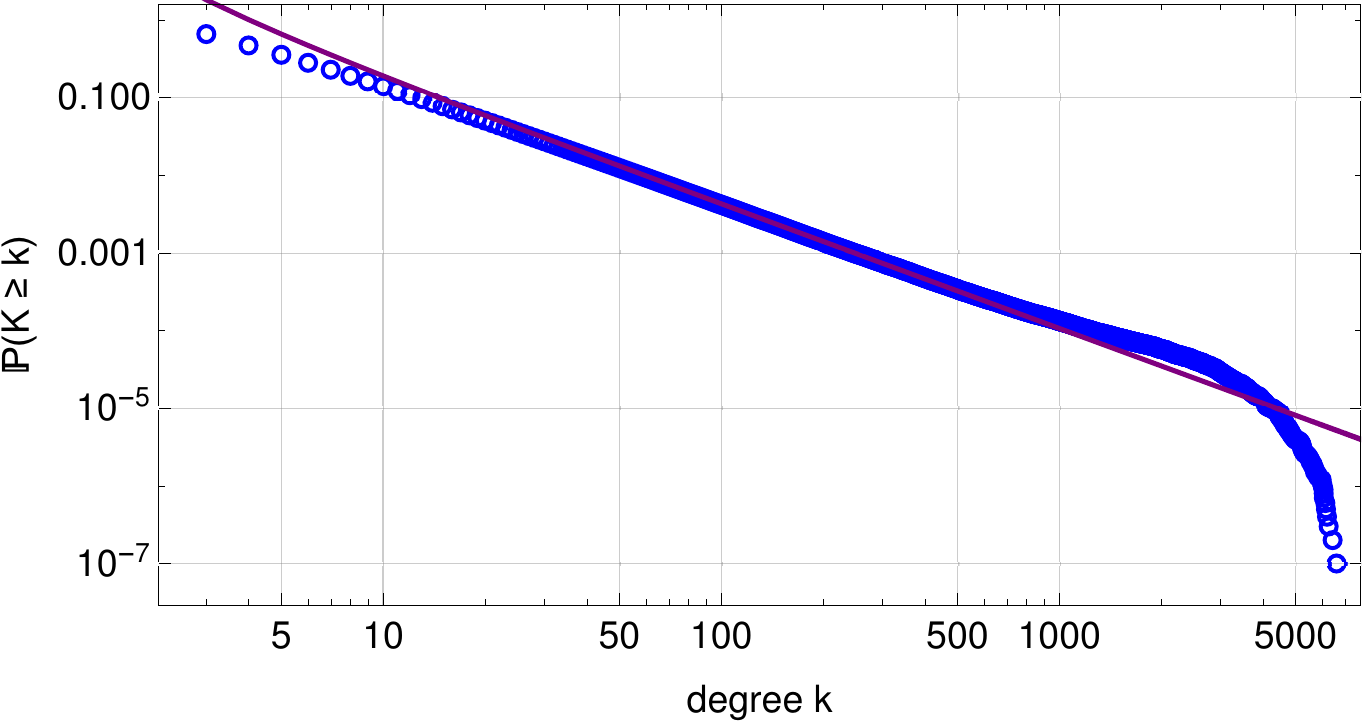}}
\caption{Degree distributions and complementary cumulative degree
  distributions for three sequences of networks. The solid lines
  represent the average of the ccdf of 100 runs of the model; dashed
  lines represent the predictions for $m = 5, \hat{m} = 3$ and $\alpha
  = 0.0$ in (a) and (b); $\alpha = 0.6$ in (c) and (d); $\alpha=1.0$
  in (e) and (f).}
\label{fig:figure4}
\end{figure}

Consider the plots in Figure~\ref{fig:figure4}. They summarize
experiments performed on three sequences of networks generated by the
algorithm in Definition~\ref{def.netmodel.netmodel} from the complete
graph with 3 nodes. In the three sequences the parameters $m = 5$ and
$\hat{m} = 3$ are fixed. However, each sequence uses a different
attachment parameter: $\alpha=0.0$, $\alpha=0.6$, and $\alpha=1.0$,
respectively.  The plots in Figure~\ref{fig:figure4} summarize the
degree distribution and the complementary cumulative degree
distribution for each one of the three sequences of networks. The main
observation is that the simulated distributions approach the
theoretical limits, a result that follows from Theorem~\ref{thm2} and
Corollary~\ref{cor.indegree.behavior}.

\section{Results}
\label{section5}

This section showcases an application of Theorem~\ref{thm1} and the EM
Algorithm (i.e., Algorithm~\ref{alg:em}). It validates the proposed
approach to estimate the contribution of the attachment mechanisms in
an empirical citation network.

The High Energy Physics Theory (HEP-Th) citation network, a publicly
available dataset compiled by ArXiv for the KDD Cup 2003
competition~\cite{Leskovec2007}, is used. This dataset covers a
network of 27770 papers and 352807 citations among them.  Each paper
has a unique identifier and is annotated with a time-stamp
corresponding to its publication date.

The seed network $G_0 = (V_0, E_0)$ has as nodes $V_0$ the papers
published in February 1992 and the papers cited by them.  There is a
directed edge $(u,v) \in E_0$ from paper $u$ to paper $v$ if and only
if paper $u$ cites paper $v$, for any $u, v \in V_0$.  In total, the
seed network $G_0$ has 4 nodes and 2 directed edges.  The sequence
$(G_t)_{t\leq 24284}$ is built by sorting the papers in the HEP-Th
network not present in $V_0$ in ascending order by time-stamps. That
is, the first published paper after February 1992 comes first and
identifies the time $t=1$, the second one identifies the time $t=2$,
and so on. There is a total of $T=24284$ papers considered in this
sequence. Starting from the seed network $G_0$, at each time step $t >
0$, the network $G_t = (V_t, E_t)$ is constructed as follows:

\begin{itemize}
  \item the set of nodes $V_t$ contains all nodes in $V_{t-1}$, and
  has a new node $u_t$ representing the paper added at time $t$ and a
  new node for each paper not present in $V_{t-1}$ that is cited by
  $u_t$; and

\item the set of edges $E_t$ contains all edges in $E_{t-1}$ and adds
  a new edge $(u_t,v)$ for each paper $v \in V_t$ cited by $u_t$.
\end{itemize}
This process results in the sequence $(G_t)_{t\leq 24284}$ of
networks, with $|V_{24284}| = 27770$ and $|E_{24284}| = 352807$,
representing the growth process of the HEP-Th citation network from
February 1992 to April 2003.

\begin{figure}[htbp]
    \centering
    \subfloat[][]{\includegraphics[width=170pt]{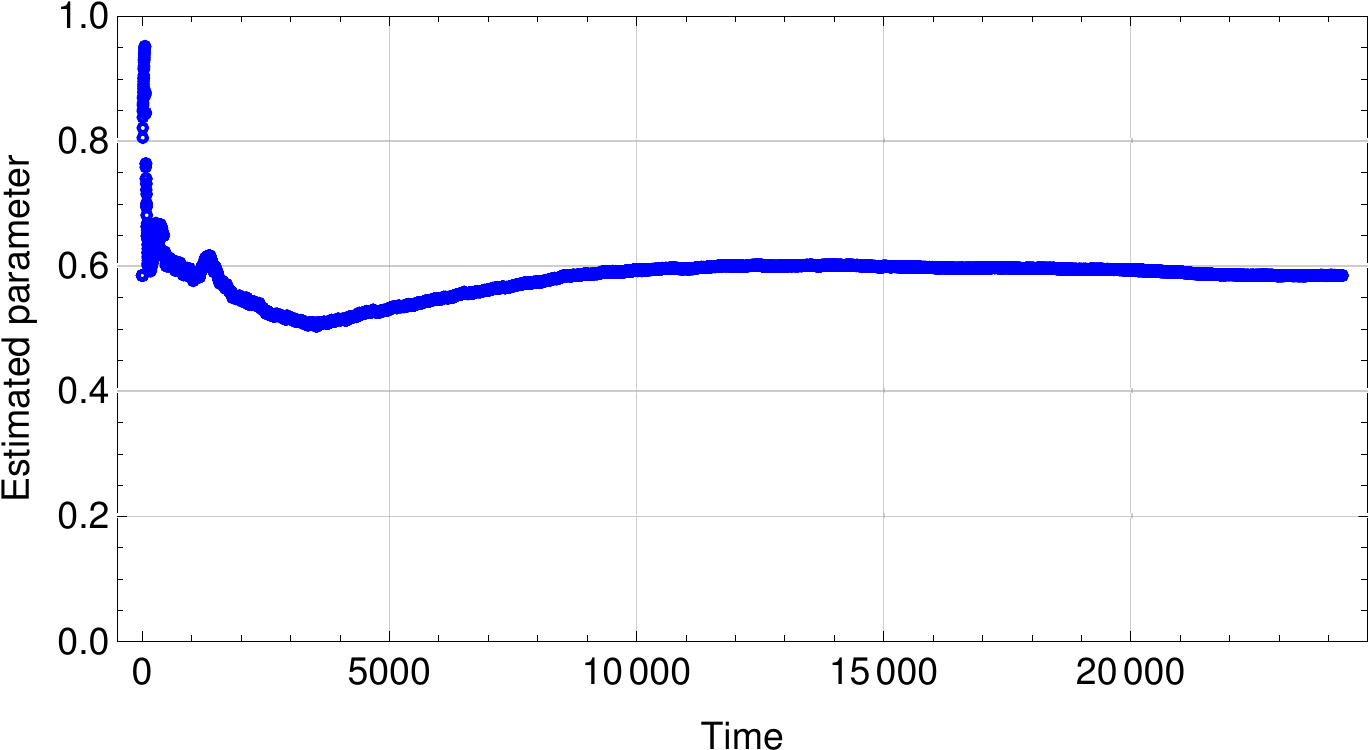}}~
    \subfloat[][]{\includegraphics[width=170pt]{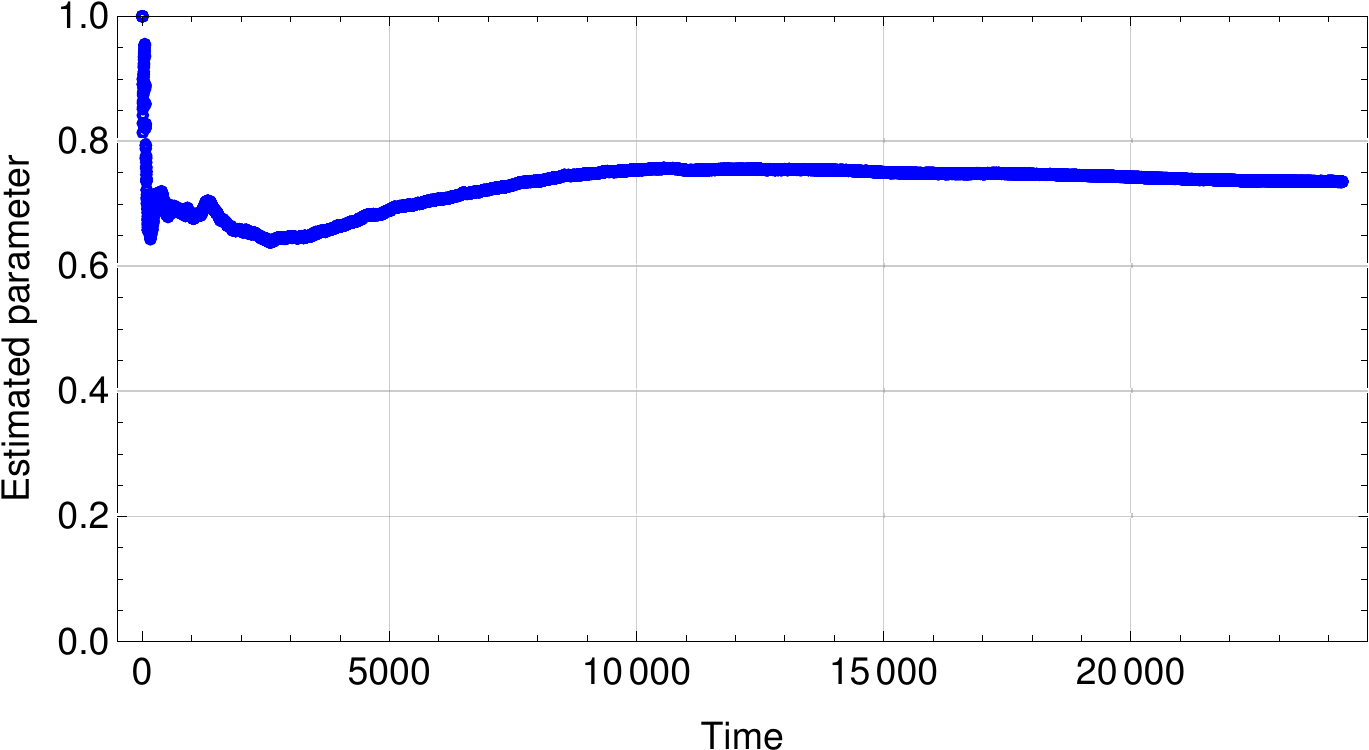}}
\caption{Evolution of the estimated parameter for the HepTh network: (a) applying Theorem~\ref{thm1} and (b) applying Algortihm~\ref{alg:em}.}
\label{fig:HepTh}
\end{figure}

The incidence proportion of the preferential attachment mechanism is
estimated by maximizing the $G_t$-likelihood at each time $t$, both by
using Theorem~\ref{thm1} and the EM Algorithm (i.e.,
Algorithm~\ref{alg:em}). The estimation of the attachment mechanisms
for $(G_t)_{t\leq 24284}$ is presented in~\Fref{fig:HepTh}, where dots
represent the estimated values for the $G_t$-likelihood function at
each time $t$. By incrementally using Theorem~\ref{thm1}, the
estimated parameter is calculated to be $\hat{\alpha}_1 = 0.59$. By
incrementally using Algorithm~\ref{alg:em}, the parameter is
calculated to be $\hat{\alpha}_2 = 0.74$. With both methods the
estimated value becomes stable around time $t=10000$, despite the fact
that this value is different for each method. One conjecture that may
explain the difference between $\hat{\alpha}_1$ and $\hat{\alpha}_2$
is related to the way 0-in-degree nodes are accounted for in each
case. On the one hand, 0-in-degree nodes contribute roots to the
$G_t$-likelihood function (Definition~\ref{def.maxlike.likelihood})
used in Theorem~\ref{thm1}.  On the other hand, the E-step of the
EM-algorithm ignores 0-in-degree nodes; thus, these nodes do not
ultimately contribute to the average computed by the
algorithm. Nevertheless, it is important to note that 0-in-degree
nodes, in general, could contribute useful information for parameter
estimation.  In the case of the HEP-Th citation network the
0-in-degree nodes identify papers that are not cited. It remains as
part of the future work to study extensions of the EM Algorithm that
consider 0-in-degree nodes.

\begin{figure}[htbp]
    \centering
    \subfloat[][]{\includegraphics[width=170pt]{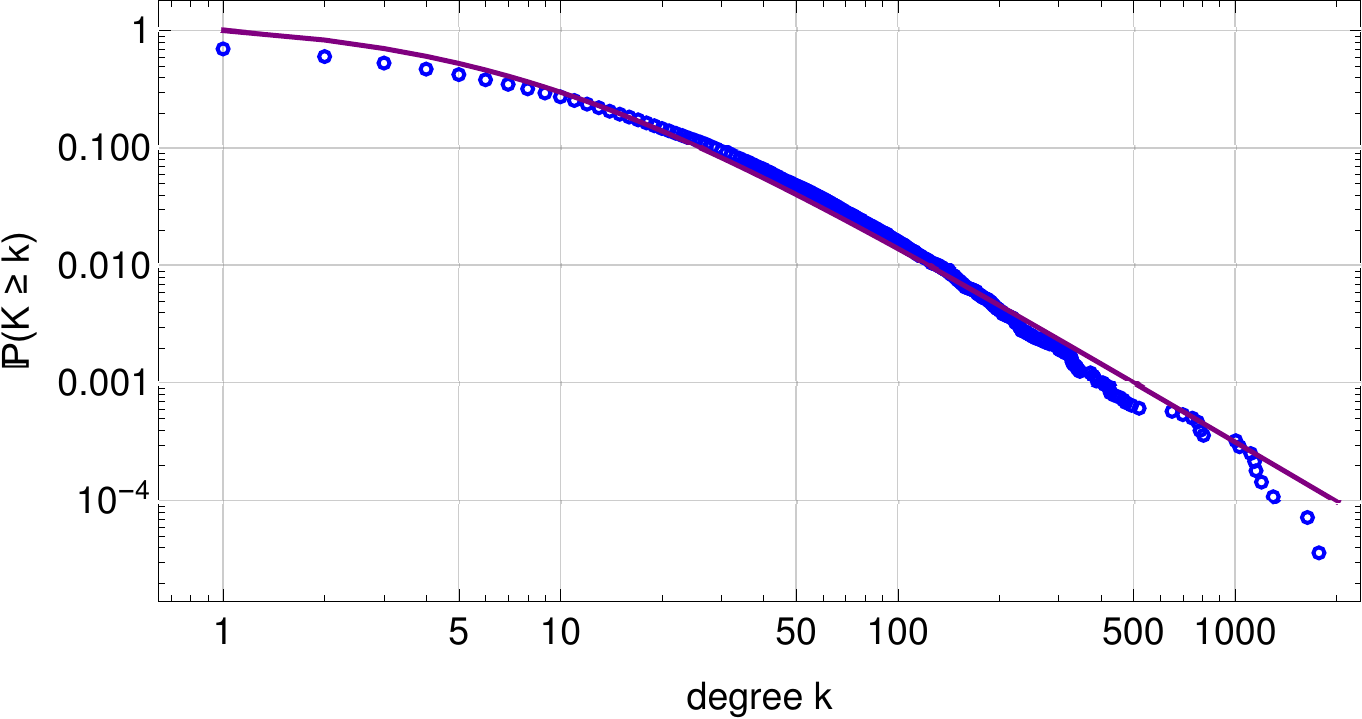}}~
    \subfloat[][]{\includegraphics[width=170pt]{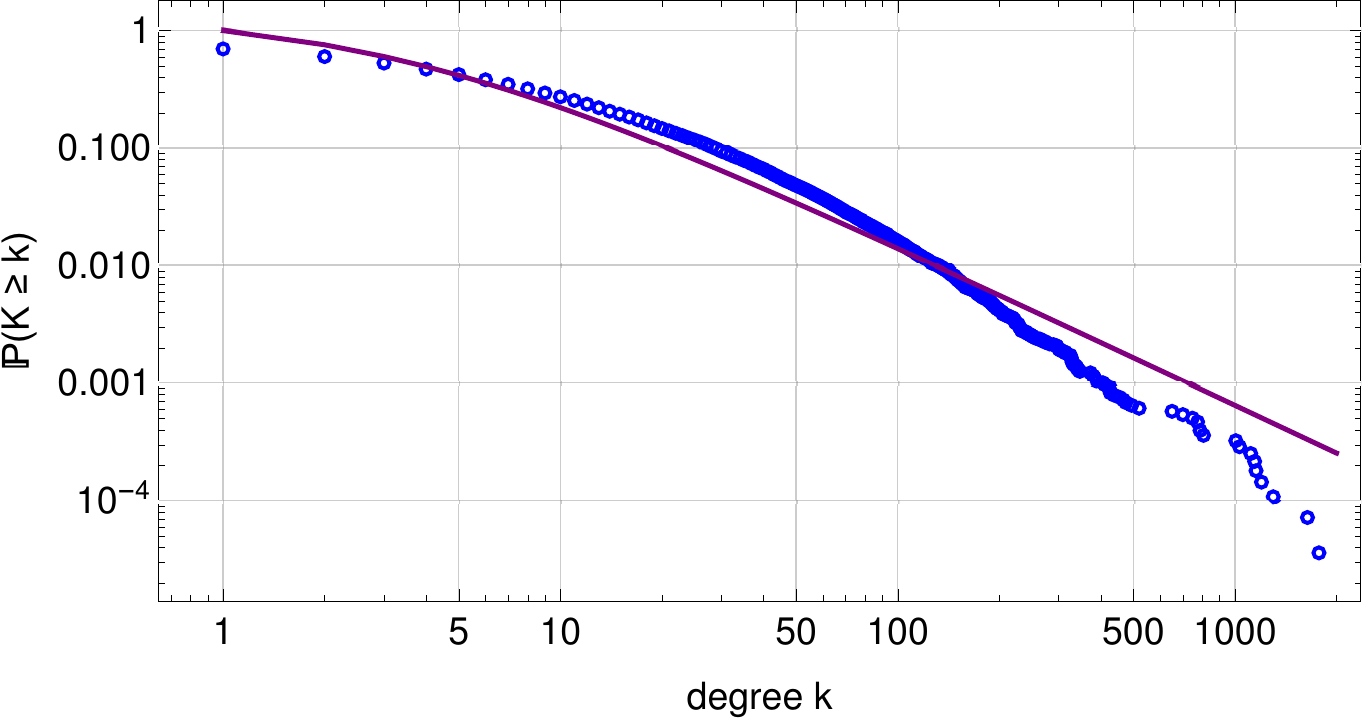}}
\caption{Complementary cumulative degree distribution for the HEP-Th
  network. Solid lines represent the predictions and dots represent
  the actual values for the empirical complementary degree
  distributions. (a) $m = 12$, $\hat{m}=0$, and $\hat{\alpha}_1=
  0.59$; (b) $m = 12$, $\hat{m}=0$, and $\hat{\alpha}_2= 0.74$.}
\label{fig:RHepTh}
\end{figure}

Based on the estimated parameters $\hat{\alpha}_1$ and
$\hat{\alpha}_2$, Theorem~\ref{thm2} and
Corollary~\ref{cor.indegree.behavior} are applied to find the
theoretical in-degree distribution. The value of the parameter $m$ is
estimated using the empirical degree distribution and the estimated
parameters; in this case study, $m$ is found to be 12. Note that,
because of the form the network is built, the parameter $\hat{m}$ is 0
(i.e., the networks do not respond to incoming nodes).
\Fref{fig:RHepTh} illustrates the relationship between the theoretical
and empirical complementary cumulative degree distributions on
$G_{24284}$.

The estimated values have a better fit when the attachment estimate is
$\hat{\alpha}_1 = 0.59$. In this case, the fit can serve as a witness
of the fact that mixed attachment models can be used to recreate the
behavior experimental networks. The case of the complementary
cumulative degree distribution when using the estimate $\hat{\alpha}_2
= 0.74$ may suggest that the actual contribution of the 0-in-degree
nodes, ignored by the EM Algorithm, weights negatively against the
estimation process.

\section{Conclusion and Future Work}
\label{section6}

Preferential attachment models explain the formation of power-laws in
the tail of degree distribution. Such models capture the evolution of
the number of connections of a small --~yet significant~-- number of
nodes with extremely large degrees. However, preferential attachment
alone falls short in describing the behavior of the large majority of
nodes with smaller degrees.

To overcome this limitation, mixed attachment models contemplate how
degree distributions may result from a combination of multiple
mechanisms. Our work is novel for it presented conditions guaranteeing
that the prevalence estimate of preferential and random attachment
mechanisms represents a local maximum of the likelihood function. We
used the expectation maximization algorithm to find the
maximum-likelihood estimate of the contribution of the two
mechanisms. Our results showed that if the algorithm is applied
without satisfying the proposed conditions, then the estimate fails to
converge to a stationary value.

Finally, we applied the proposed approach to estimate the prevalence
of random and preferential attachment mechanism in citation networks
of academic papers. The estimate is evaluated by comparing the
empirical degree distribution to the theoretical distribution
evaluated at the estimated parameter. The results showed that mixed
attachment models are able to recreate the behavior of nodes with both
small and large degrees.

Future work on extending the proposed model to include new attachment
and response mechanisms that can update edges and even generate
clustering should be pursued. For instance, considering link creation
rates based on out-degree distributions (e.g., as proposed
in~\cite{krapivsky2001}) and internal rewiring (e.g., as proposed
in~\cite{krapivsky2001b}), should be considered.  Furthermore, the
analysis of the likelihood functions, and the in- and out-degree
distributions of the extended models should also be investigated. New
applications to other empirical networks should also be considered,
taking into account the rich experience already available with mixed
attachment networks. Finally, extensions of the EM Algorithm that take
into account 0-in-degree nodes must be studied in order to apply these
techniques to, e.g., tree-like empirical networks.

\paragraph{Acknowledgments.}
Thanks are due to D. Ruiz and F. Amaya for valuable feedback and
fruitful discussions. The three authors were in part supported by the
Ministry of Information Technologies and Telecommunications of
Colombia (MinTIC), and the Colombian Administrative Department of
Science, Technology and Innovation (COLCIENCIAS) through the Center of
Excellence and Appropriation in Big Data and Data Analytics (CAOBA).

\bibliographystyle{abbrv}
\bibliography{biblio}

\end{document}